\documentclass{article}

\usepackage{arxiv}

\usepackage[utf8]{inputenc} 
\usepackage[T1]{fontenc}    
\usepackage{hyperref}       
\usepackage{url}            
\usepackage{booktabs}       
\usepackage{amsbsy,amsfonts,amsmath,amssymb}
\usepackage{algorithm}
\usepackage{algpseudocode}
\usepackage{tikz}
\usepackage{nicefrac}       
\usepackage{microtype}      
\usepackage{lipsum}		
\usepackage{graphicx}
\usepackage{natbib}
\usepackage{doi}
\usepackage{placeins} 
\usepackage{setspace} 
\usepackage{bm}

 \usepackage{array,multirow,graphicx}
 \usepackage{float}
 \usepackage{arydshln}

\usepackage{natbib}
 \bibpunct[, ]{(}{)}{,}{a}{}{,}%
\usepackage{caption}
\usepackage[labelfont=sf]{subcaption}
\usepackage{tikz}
\usetikzlibrary{positioning,shapes,decorations.pathmorphing}
\usetikzlibrary{patterns}
\usepackage{pgfplots}
\pgfplotsset{compat=1.8}
\usepgfplotslibrary{statistics}
\usepackage{pgfplots}
    \usepgfplotslibrary{
        groupplots,
        units,
    }
    \pgfplotsset{compat=1.3}
\usetikzlibrary{arrows.meta}
\usetikzlibrary{shapes.geometric}
\usepackage{setspace}
\pgfplotsset{compat=1.12}

\title{Enhancing Kernel Search with Pattern Recognition: the Single-Source Capacitated Facility Location Problem}

\author{ \href{https://orcid.org/0000-0002-9340-5407}{\includegraphics[scale=0.06]{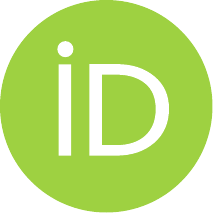}\hspace{1mm}Hannah ~Bakker}\\
	Institute of Operations Research\\
	Karlsruhe Institute of Technology\\
	Karlsruhe, Germany \\
	\texttt{hannah.bakker@kit.edu} \\
	\And
	\href{https://orcid.org/0000-0002-7810-4112}{\includegraphics[scale=0.06]{orcid.pdf}\hspace{1mm}Gianfranco ~Guastaroba} \\
	Department of Economics and Management\\
	University of Brescia\\
	Brescia, Italy \\
	\texttt{gianfranco.guastaroba@unibs.it} \\
	  \And
	\href{https://orcid.org/0000-0002-8339-0117}{\includegraphics[scale=0.06]{orcid.pdf}\hspace{1mm}Stefan ~Nickel} \\
	Institute of Operations Research\\
	Karlsruhe Institute of Technology\\
	Karlsruhe, Germany \\
	\texttt{stefan.nickel@kit.edu} \\
    \And
	\href{https://orcid.org/0000-0002-8893-5227}{\includegraphics[scale=0.06]{orcid.pdf}\hspace{1mm}M. Grazia ~Speranza} \\
	Department of Economics and Management\\
	University of Brescia\\
	Brescia, Italy \\
	\texttt{grazia.speranza@unibs.it} \\
}

\date{}

\renewcommand{\headeright}{ }
\renewcommand{\undertitle}{ }
\renewcommand{\shorttitle}{PaKS: Pattern-based Kernel Search}

\hypersetup{
pdftitle={Enhancing Kernel Search with Pattern Recognition: the Single-Source Capacitated Facility Location Problem},
pdfsubject={q-bio.NC, q-bio.QM},
pdfauthor={Hannah ~Bakker, Gianfranco ~Guastaroba, Stefan ~Nickel, M. Grazia ~Speranza},
pdfkeywords={Pattern recognition, Biclustering, Kernel Search, Single-source capacitated location problems, Heuristic framework},
}

\newcommand{\subjto}{\ensuremath{\text{subject to}}}
\newcommand{\prob}{{\ttfamily{SSCFLP}}}
\newcommand{\model}{{\ttfamily{SSCFL model}}}
\newcommand{\alg}{{\ttfamily{PaKS}}}
\newcommand{\ks}{{\ttfamily{KS14}}}
\newcommand{\cpx}{{\ttfamily{CPLEX}}}
\newcommand{\hils}{{\ttfamily{HILS}}}
\newcommand{\mscflp}{{\ttfamily{MSCFLP}}}
\newcommand{\medcirc}{\raisebox{0.25pt}{\scalebox{0.8}{$\circ$}}}

\begin{document}
\maketitle

\begin{abstract}
We introduce Pattern-based Kernel Search (\alg), a two-phase matheuristic for the solution of the Single-Source Capacitated Facility Location Problem (\prob). In the first phase, \alg\ employs a pattern recognition technique to identify an implicit spatial separation of potential locations and customers into subsets, called regions, within which location and assignment decisions are strongly interdependent. In the second phase, \alg\ employs an enhanced Kernel Search (KS) heuristic that leverages the interdependencies among the decision variables identified in the first phase. On a set of 112 benchmark instances, consisting of up to 1,000 locations and 1,000 customers, computational results show that \alg\ consistently outperforms both a standard KS implementation and the current state-of-the-art heuristic for solving the \prob, as well as \cpx\ when run with a time limit. For these instances, \alg\ achieved an average gap compared to the best known solution of 0.02\%. Experimental results conducted on a large set of new very large test problems, comprising up to 2,000 locations and 2,000 customers, demonstrate that \alg\ outperforms both the standard KS heuristic and \cpx\ in terms of quality of the solution found, finding the largest number of best solutions, and achieving the smallest average gap.

\end{abstract}

\keywords{Pattern recognition, Biclustering, Kernel Search, Single-source capacitated location problems, Heuristic framework}

\section{Introduction}\label{sec:Intro}

Facility location problems call for determining an optimal location of one
or multiple facilities in order to serve the demand of a given set of customers. It is among the most investigated classes of combinatorial optimization problems, and attracts an ever-growing interest among operational researchers and practitioners. The book edited by \citet{laporte2020location} provides a comprehensive coverage of both basic and advanced problems, formulations, and algorithms devised for this active research area.

The \textit{Single-Source Capacitated Facility Location Problem} (\prob) is among the prototypical optimization problems in the class of discrete location problems, for which several exact and heuristic algorithms have been proposed. In the \prob, we are given a set $J$ of demand nodes, hereafter called customers. Each customer $j \in J$ is associated with a demand $d_j$ to be served. We are also given a set $I$ of potential locations where a facility can be open. For the sake of brevity, henceforth we refer to a potential location $i \in I$ simply as location $i$. Each location $i \in I$ is associated with a fixed opening cost $f_i$ and a capacity $q_i$. The cost of assigning one unit of demand of customer $j$ to location $i$ is denoted by $c_{ij}$. The \prob\ calls for selecting an optimal subset of locations to open from $I$, along with an optimal assignment of all customers in $J$ to the open locations, so that the sum of the opening and assignment costs is minimized. Each customer demand must be entirely supplied by one location, and each open location cannot supply an amount of demand greater than its capacity. 

The single-sourcing restriction makes the \prob\ particularly challenging to solve \citep[e.g., see][]{klose2005facility}. In fact, for a given set of open locations, the problem of assigning customers to them is a special case of the Generalized Assignment Problem, which is $\mathcal{NP}$-hard by itself \citep[see][]{nauss2003solving}. Nevertheless, the single-source assumption formalizes a number of features observed in several real-world applications \citep[e.g., see][]{guastaroba2014heuristic}.

Let $y_i \in \{0,1\}$ be a binary variable taking value 1 if location $i \in I$ is open, and 0 otherwise. Let $x_{ij} \in \{0, 1\}$ be a binary variable that takes value 1 if the demand of customer $j \in J$ is assigned to location $i \in I$, and 0 otherwise. Then, the \prob\ can be cast as the following pure Binary Integer linear Program (BIP):
\begin{align}
& [\text{\model}] & \nonumber\\
\min \quad & \sum_{i \in I} f_i y_i + \sum_{j \in J} d_j \sum_{i \in I}  c_{ij}  x_{ij} \label{SSCFL-obj} \\
\subjto{} \quad & \sum_{i \in I} x_{ij} = 1 & (j \in J) \label{SSCFL-assignment}\\
           &  \sum_{j \in J} d_j x_{ij} \leq q_i y_i  & (i \in I) \label{SSCFL-demand}\\
           &  x_{ij} \leq y_i  & (i \in I;  j \in J) \label{SSCFL-VI}\\
           &  y_{i} \in \{0,1\} & (i \in I) \label{SSCFL-yBin}\\
           &  x_{ij} \in \{0,1\} & (i \in I;  j \in J).\label{SSCFL-xBin}
\end{align}
Objective function~\eqref{SSCFL-obj} calls for the minimization of the sum of two components: the total fixed cost of opening the locations and the total demand assignment cost, respectively. Constraints~\eqref{SSCFL-assignment} impose that the demand of each customer is completely served and, along with constraints~\eqref{SSCFL-xBin}, entirely allocated to a single location. Constraints~\eqref{SSCFL-demand} ensure that, on the one hand, the capacity of each open location is not violated, and, on the other hand, prevent customers to be allocated to non-open locations. Constraints~\eqref{SSCFL-yBin} and~\eqref{SSCFL-xBin}
define the binary domain of the decision variables. Note that inequalities~\eqref{SSCFL-VI}, being redundant in the \model, are well-known to yield a much tighter Linear Programming (LP) relaxation than the equivalent formulation without them \citep[see, among others,][]{yang2012cut}. Nonetheless, the number of inequalities~\eqref{SSCFL-VI} increases very rapidly as the size of the instance grows. \citet{guastaroba2014heuristic} point out that, as a consequence, it might be challenging to solve, in reasonable computing times, the LP relaxation of large \prob\ instances with an off-the-shelf solver if all inequalities~\eqref{SSCFL-VI} are included upfront. Finally, without loss of generality, we assume that $d_j > 0$ for all $j \in J$, $f_i > 0$ and $q_i > 0$ for all $i \in I$, $c_{ij} \geq 0$ for all $i \in I$ and $j \in J$.

It is beyond our scope to provide an extensive overview of the literature on the \prob. Thus, we limit the review to the papers that are the most recent and closely related to our research. The reader interested in a more general overview is referred to~\citet{fernandez2019fixed} and to the references cited therein, as well as in the papers cited below.

The \textit{focus of this paper} is on developing a two-phase matheuristic for solving the \prob\ that combines a pattern recognition technique with a Kernel Search (KS) heuristic. We call the resulting algorithm the \textit{Pattern-based Kernel Search} (\alg). In general terms, pattern recognition refers to a class of data analysis methods that use algorithms to automatically identify patterns and regularities in data. These methods have been applied to a wide variety of data, from text and images to sounds. Potentially, they can also be applied to a given optimization problem in order, for example, to identify structural patterns characterizing an optimal (or near-optimal) solution. These patterns could then be exploited by specifically designed solution algorithms. Despite these potential advantages, the application of such techniques to determine structural patterns in optimal solutions is essentially unexplored. To our knowledge, \citet{Bakker2024} is the only author that, in the domain of facility location problems, employs pattern recognition techniques to identify the implicit spatial structure of an instance of a Multi-Source Capacitated Facility Location Problem (\mscflp). The author applies a technique, called spectral biclustering, to group the sets of locations and customers into subsets that can be interpreted as regions, each showing strong internal interdependencies. The purpose is to determine subsets such that across high-quality solutions customers are generally allocated to locations within the same region, while the assigned locations may vary from one solution to another. As a consequence, this implies that the associated decision variables interdepend more strongly within a region than across different regions. The interpretation of a subset as a region is intuitive given the spatial dimension underlying the location problem.

 KS is an integer programming based heuristic framework that has been shown to produce high-quality solutions for several classes of combinatorial optimization problems, including knapsack problems \citep[][]{angelelli2010kernel, lamanna2022two}, routing problems \citep[][]{archetti2021kernel, gobbi2023hybridizing}, feature selection problems in support vector machine models \citep[][]{labbe2019mixed}, lot sizing problems \citep[][]{carvalho2018kernel}, supply chain optimization problems \citep[][]{zhang2019novel}, maximum flow problems \citep[][]{carrabs2025hybridizing}, and general Mixed Integer linear Programs (MIPs) \citep[][]{guastaroba2017adaptive}. The reader interested in a general description of the approach is referred to the tutorial provided in~\citet{maniezzo2021kernel}.

KS has been particularly successful in solving discrete location problems \citep[see][]{guastaroba2012kernel, santos2020kernel, filippi2021kernel}. Related to the present research, \citet{guastaroba2014heuristic} implement a conventional KS to solve the \prob. Henceforth, this heuristic is referred to as the \ks. Later, \citet{tran2017hypergraph} solve the \prob\ by using a large neighborhood search heuristic that exploits specifically designed neighborhood structures. \citet{caserta2020general} introduce a matheuristic algorithm based on the corridor method to solve a variety of capacitated location problems, including the \prob. To our knowledge, the current state-of-the-art heuristic for solving the \prob\ is the Hybrid Iterated Local Search (\hils) devised in \citet{de2024hybrid}.

\noindent \textbf{Contributions of the paper.} The key contributions of the present research can be summarized as follows. 

\begin{itemize}
    \item[\checkmark] From an algorithmic perspective, we introduce \alg, a two-phase matheuristic for the solution of the \prob. In the first phase, \alg\ employs a pattern recognition technique to gather information on the importance and interdependence among decision variables, grouping locations and customers into regions. \alg\ leverages this information to determine the sequence of restricted models that are solved in the second phase. Thereby, it expands the conventional KS algorithm, which determines this sequence solely based on the solution of the LP relaxation, to include a more extensive information gathering phase that combines techniques from mathematical programming with approaches from data science. As a result, the size of the restricted models can be substantially reduced, which proved to be particularly beneficial while solving large-scale instances.
    
    \item[\checkmark] Regarding the solution of the \prob, this paper presents a robust, reliable, and effective heuristic for solving large and very large-scale instances. Extensive computational experiments conducted on benchmark instances, consisting of up to 1,000 locations and 1,000 customers, show that \alg\ outperforms the above-mentioned state-of-the-art heuristic developed in  \citet{de2024hybrid}, as well as a re-implementation of \ks\ and a general-purpose MIP solver --namely, \cpx-- when the computing time is limited. Additional computational experiments carried out on new and very large-scale instances, which comprise up to 2,000 locations and 2,000 customers, point out that \alg\ significantly outperforms both the re-implementation of \ks\ and \cpx.
\end{itemize}
    
\noindent \textbf{Structure of the paper.} The paper is organized as follows. \alg\ is introduced and described in detail in Section~\ref{sec:Method}. The results of an extensive experimental analysis are summarized and discussed in Section~\ref{sec:Experiments}. Finally, some
concluding remarks are outlined in Section~\ref{sec:Conclusions}.

\section{The Pattern-based Kernel Search}\label{sec:Method}

In this section, we describe \alg, the matheuristic developed for the solution of the \prob. We introduce some basic definitions and notation in Section~\ref{sec:Defs}. In Section~\ref{sec:PKS}, we illustrate the \alg\ providing an intuitive description of the approach. Then, we detail the two main phases composing the \alg\ in Sections~\ref{sec:Phase1} and \ref{sec:Phase2}, respectively, emphasizing the foremost changes compared to a standard KS implementation as the \ks.

\subsection{Preliminary definitions and basic notation}\label{sec:Defs}

Let $\mathcal{Y}$ denote the set of all binary variables $y_i$, with $i \in I$. Similarly, let $\mathcal{X}$ be the set of all $x_{ij}$, with $i \in I$ and $j \in J$. We denote as $U_y \subsetneqq \mathcal{Y}$ and as $U_x \subsetneqq \mathcal{X}$ a subset of the $y$ and of the $x$ binary variables, respectively. We denote as $BIP(U_y \cup U_x)$ the \textit{restricted} BIP in which the $y$ variables in $\mathcal{Y} \setminus U_y$ and the $x$ variables in $\mathcal{X} \setminus U_x$ are all fixed to zero. Additionally, let $z^{H}$ be an upper bound to the value of an optimal solution to the \textit{original problem}, denoted by $BIP(\mathcal{Y} \cup \mathcal{X})$. In the following, we say that a binary variable is \textit{promising} if its value is presumably one in an optimal solution to the original problem. Further, let $LP(\mathcal{Y} \cup \mathcal{X})$ be the LP relaxation of the original problem, and $LP(U_y \cup U_x)$ denote the LP relaxation restricted to the variables in $U_y$ and $U_x$ --i.e., the LP relaxation of $BIP(U_y \cup U_x)$. 

In general terms, the KS algorithmic framework solves, by using a general-purpose MIP solver, a sequence of restricted BIPs, each defined by a \textit{kernel} $K$ and a, possibly empty, \textit{bucket} $B^r$ with the following properties: (i) an optimal solution to $BIP(K \cup B^r)$ is an optimal or high-quality solution to the original problem; and (ii) $BIP(K \cup B^r)$ can be solved in much shorter computing times than the original problem. In \alg, kernel $K$ comprises a subset of $\mathcal{Y}$ and $\mathcal{X}$. Formally, the kernel is defined as $K = K_y \cup K_x$, where $K_y \subsetneqq \mathcal{Y}$ and $K_x \subsetneqq \mathcal{X}$. Similarly, a generic bucket is denoted as $B^r = B^r_y \cup B^r_x$, where $B_y \subseteq \mathcal{Y} \backslash K_y$ and $B_x \subseteq \mathcal{X} \backslash K_x$. To be effective, kernel $K$ should contain most of --ideally all-- the promising variables. Nevertheless, it is very difficult to identify the set of all promising variables upfront. Hence, \alg\ determines an initial composition of $K$, and then solves a sequence of BIPs, each restricted to the kernel plus, from the second BIP onward, a bucket. The solution of each BIP is used to update the composition of the kernel according to a `learn and adjust' mechanism that adds to $K$ newly discovered promising variables or remove variables that turned out to be unpromising. 

In the conventional KS, the optimal solution to $LP(\mathcal{Y} \cup \mathcal{X})$ is exploited to determine the initial composition of the kernel, as well as that of a sequence of buckets. In contrast with this standard approach, \alg\ uses information gathered from optimal solutions to $LP(\mathcal{Y} \cup \mathcal{X})$ and to a number of restricted LP relaxations. It uses pattern recognition to extract persistent allocation patterns across these solutions that group the sets of locations and customers into non-empty subsets, called \textit{regions}. Let $\mathcal{R}$ be the set of regions produced. Each region $r \in \mathcal{R}$ is defined as a tuple $(I_r, J_r)$, where $I_r \subseteq I$ denotes a subset of locations and $J_r \subseteq J$ a subset of customers.

\subsection{The algorithmic framework}\label{sec:PKS}

\begin{figure}[htbp]
\caption{General scheme of the \alg.} \label{fig:FlowChart}
\includegraphics[trim=0cm 0cm 2.0cm 0cm, clip=true, totalheight=0.29\textheight]{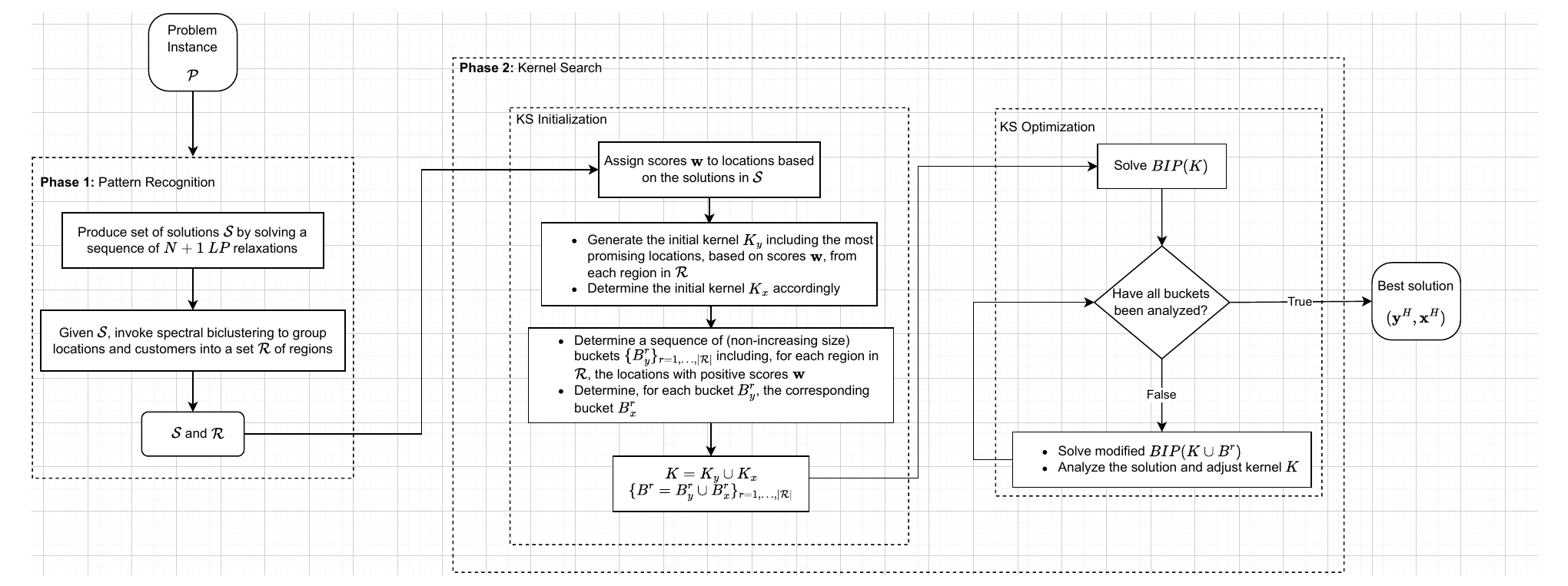}
 \end{figure}

\begin{figure}[htbp]
  \centering
  \subcaptionbox{Benchmark instance \texttt{i300\_1}.\label{fig:Inst}}{%
    \includegraphics[trim=8cm 0cm 8cm 0cm,clip,width=0.4\textwidth]{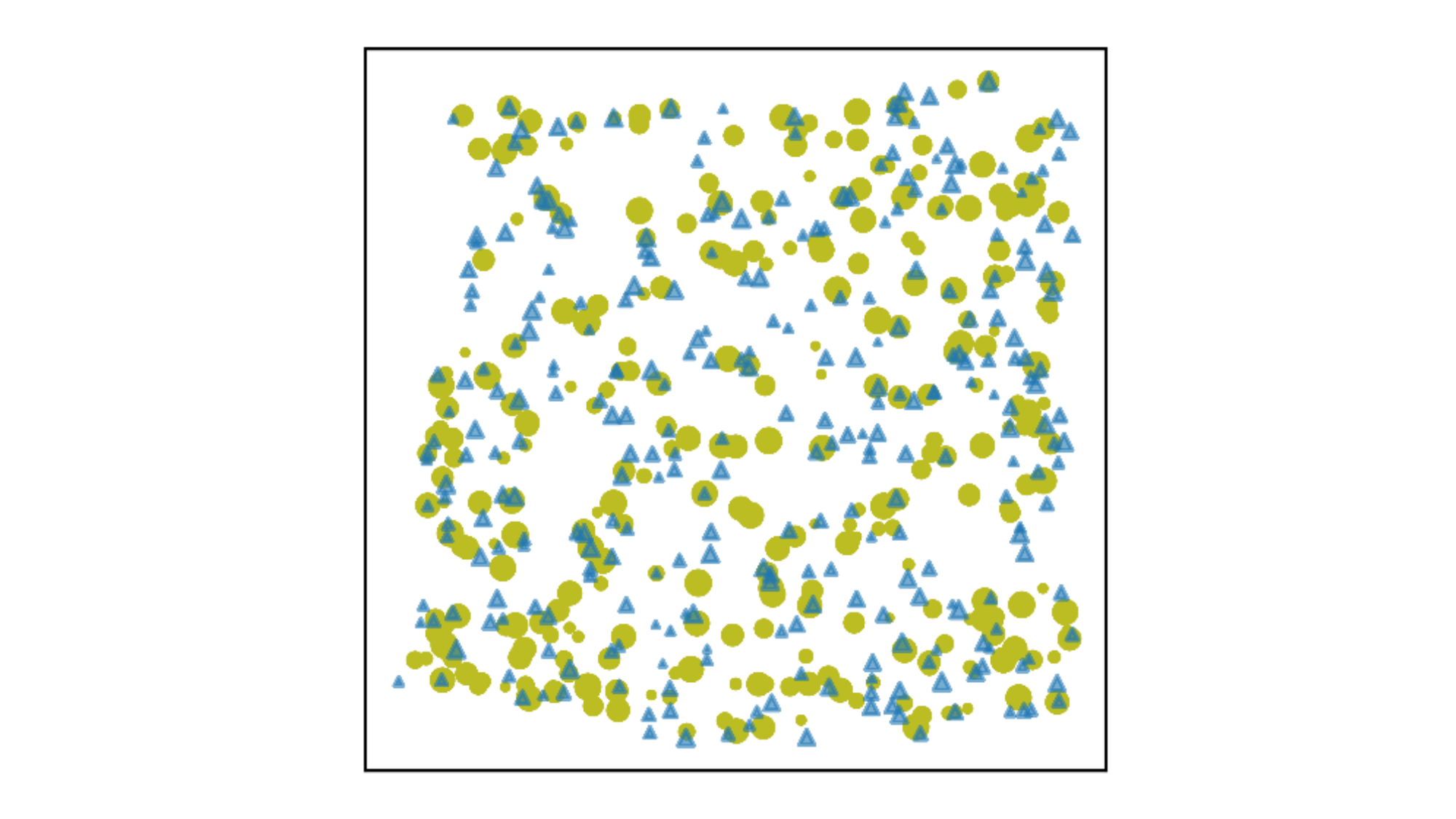}%
  }
  \hfill
  \subcaptionbox{Set $\mathcal{R}$ comprising four regions.\label{fig:Regions}}{%
    \includegraphics[trim=6cm 0cm 8cm 0cm,clip,width=0.43\textwidth]{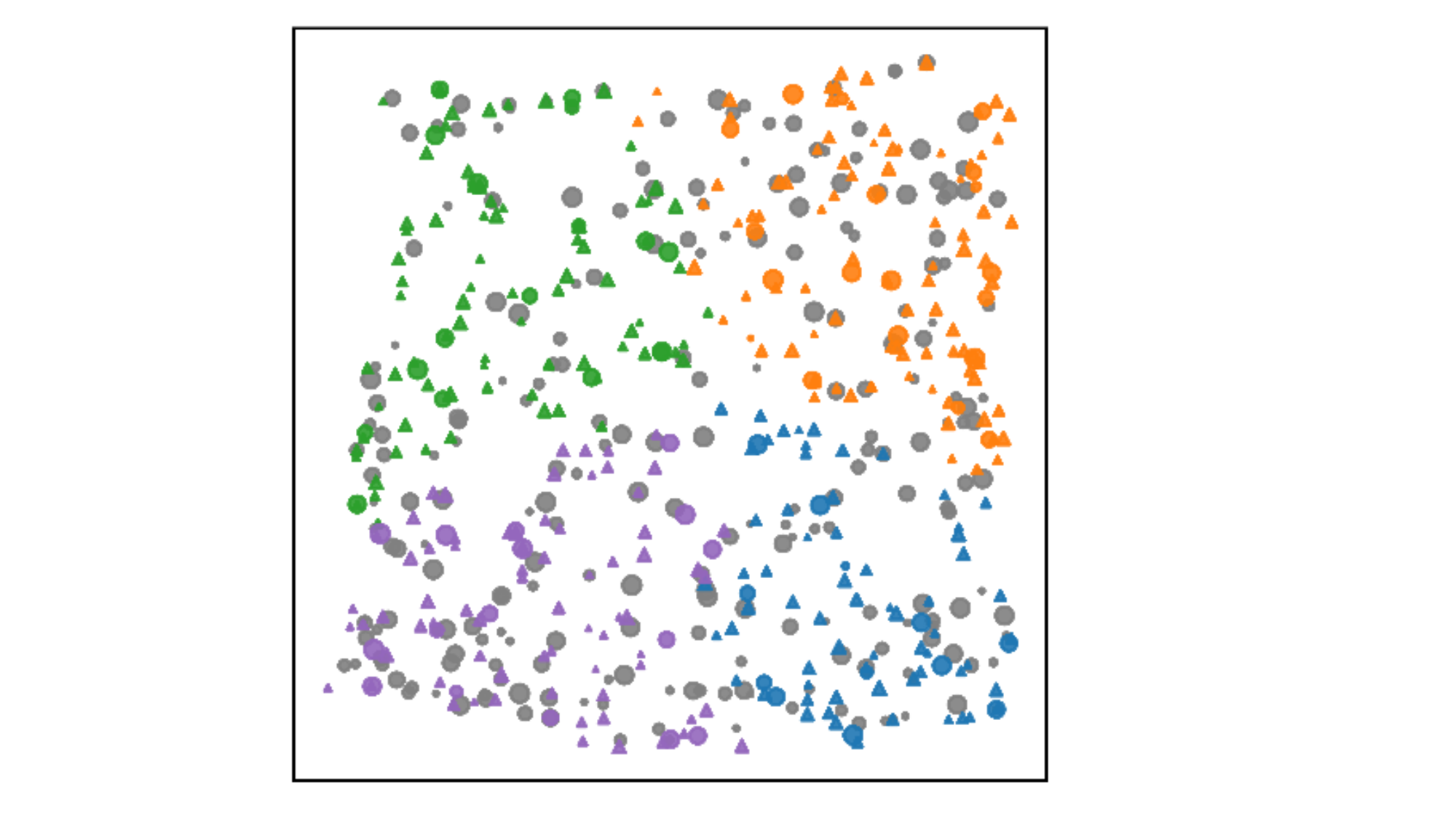}%
  }

  \caption{Instance \texttt{i300\_1} and the set of regions identified in Phase~1 of \alg. Locations are depicted as circles, whereas customers are shown as triangles.}
  \label{fig:InstRegion}
\end{figure}

This section provides a high-level description of \alg\ that follows the general algorithmic structure outlined in Figure~\ref{fig:FlowChart}.
\alg\ consists of two phases. Given an instance $\mathcal{P}$ of the \prob, \textit{Phase 1} mainly serves the purpose of gathering information on the importance of and interdependence among decision variables, and groups the locations in $I$ and the customers in $J$ into a set of regions $\mathcal{R}$. To produce such a clustering, the algorithm first solves a sequence of $N+1$ LP relaxations. The first is the LP relaxation of the original problem, that is, $LP(\mathcal{Y} \cup \mathcal{X})$. Each of the others is an LP relaxation $LP(U_y \cup U_x)$ restricted to given subsets of variables $U_y$ and $U_x$. Subsets $U_y$ are determined by randomly selecting subsets of $y$ variables chosen among those taking a nonzero value in the optimal solution found to $LP(\mathcal{Y} \cup \mathcal{X})$, and then fixing all of them to zero. For each variable belonging to a given subset $U_y$, say $y_{i^\prime}$, subset $U_x$ comprises all associated allocation variables --i.e., all $x_{i^\prime j}$, with $j \in J$. Let $\mathcal{S}$ denote the set comprising these $N+1$ solutions. Subsequently, set $\mathcal{S}$ is input to a spectral biclustering algorithm to determine the set of regions $\mathcal{R}$. Each location that is not selected in none of the former $N+1$ solutions, if any, is removed before invoking the biclustering algorithm. To this aim, in each restricted BIP solved in the remainder of the \alg, the corresponding variable $y_i$ is fixed to zero. Figure~\ref{fig:Inst} depicts benchmark instance {\texttt{i300\_1}}. Figure~\ref{fig:Regions} shows the set of regions $\mathcal{R}$ produced by the spectral biclustering algorithm. The latter method produced four regions, highlighted in green, orange, blue, and purple colors, respectively. The locations that have been removed are depicted in gray. Section~\ref{sec:Phase1} provides a more detailed description of Phase 1.

\begin{figure}[htbp]
    \centering

    \subcaptionbox{Initial kernel $K_y$ (circles) and the customers.\label{fig:InitKKS2014}}{%
        \includegraphics[trim=0cm 0cm 0cm 0cm, clip, width=0.4\textwidth]{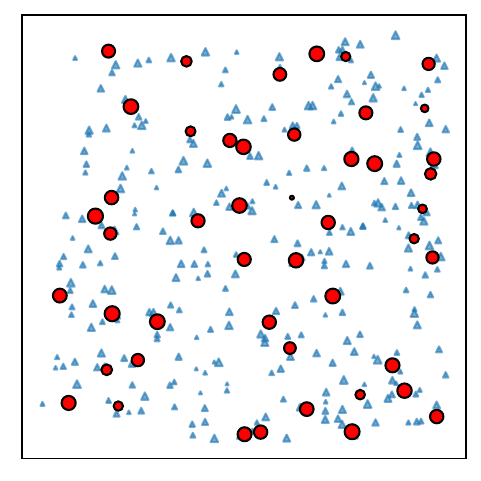}%
    }
    \hfill
    \subcaptionbox{The sequence of six buckets.\label{fig:BucketsKS2014}}{%
        \includegraphics[trim=2cm 0cm 2cm 0cm, clip, width=0.53\textwidth]{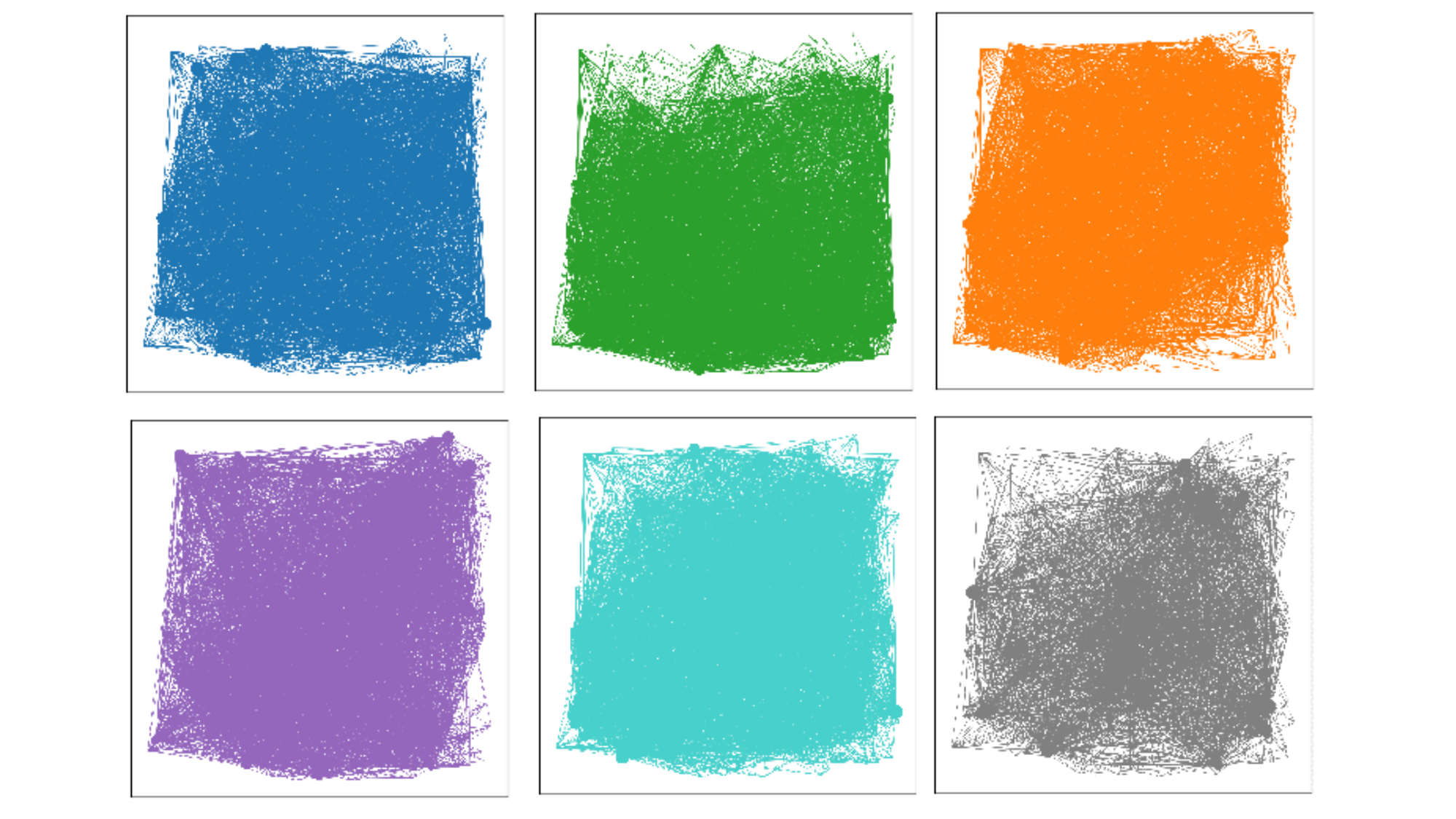}%
    }

    \caption{The initial kernel and the sequence of buckets produced by \ks\ 
    for instance \texttt{i300\_1}.}
    \label{fig:KernelBucketsKS2014}

\end{figure}

The output of Phase 1 is a set of solutions $\mathcal{S}$ and a set of regions $\mathcal{R}$. These two sets are exploited in the \textit{initialization step} of \textit{Phase 2} of \alg\ in order to determine the initial composition of the kernel and the sequence of buckets. Initially, kernel $K = K_y \cup K_x$ consists of the most promising $y$ and $x$ variables. To identify these variables, a score $w_i$, indicating the potential of location $i$, is computed based on information obtained from the solutions in $\mathcal{S}$. Then, for each region  $r \in \mathcal{R}$ the $y$ variables with the highest score are added to $K_y$. For each variable $y_i$ in $K_y$, the most promising $x$ variables are added to $K_x$ based on the information obtained from $\mathcal{S}$ and $\mathcal{R}$. Generally, the most promising $x_{ij}$ variables represent \textit{intra-regional} allocations --i.e., location $i$ and customer $j$ are located within the same region. In some circumstances, detailed below, it is beneficial to include also a small subset of the \textit{inter-regional} allocations. Then, a sequence $\{B^r_y\}_{r=1,\dots,|\mathcal{R}|}$ of buckets is determined with the $y$ variables not included in $K_y$ that have a positive score $w_i$. One bucket is created for each region. The sequence of buckets above is sorted in non-increasing order of the bucket cardinality. Next, the sequence $\{B^r_x\}_{r=1,\dots,|\mathcal{R}|}$ of buckets is determined similarly to $K_x$, adding to each bucket the most promising assignment variables.
\begin{figure}[htbp]
    \centering

    \subcaptionbox{Initial kernel $K_y$ (circles) and the customers.\label{fig:InitKPKS}}{%
        \includegraphics[trim=0cm 0cm 0cm 0cm, clip, width=0.4\textwidth]{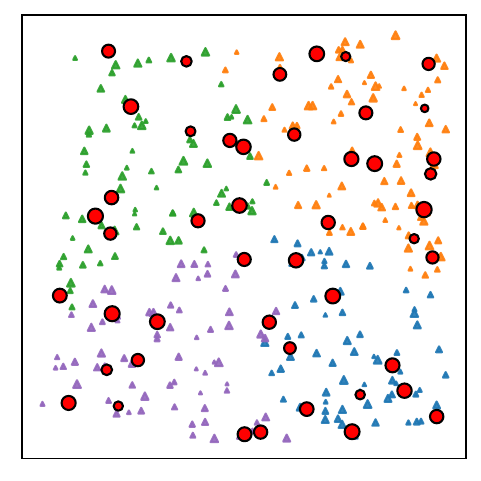}%
    }
    \hfill
    \subcaptionbox{The sequence of four buckets.\label{fig:BucketsPKS}}{%
        \includegraphics[trim=6.5cm 1cm 8cm 0cm, clip, width=0.43\textwidth]{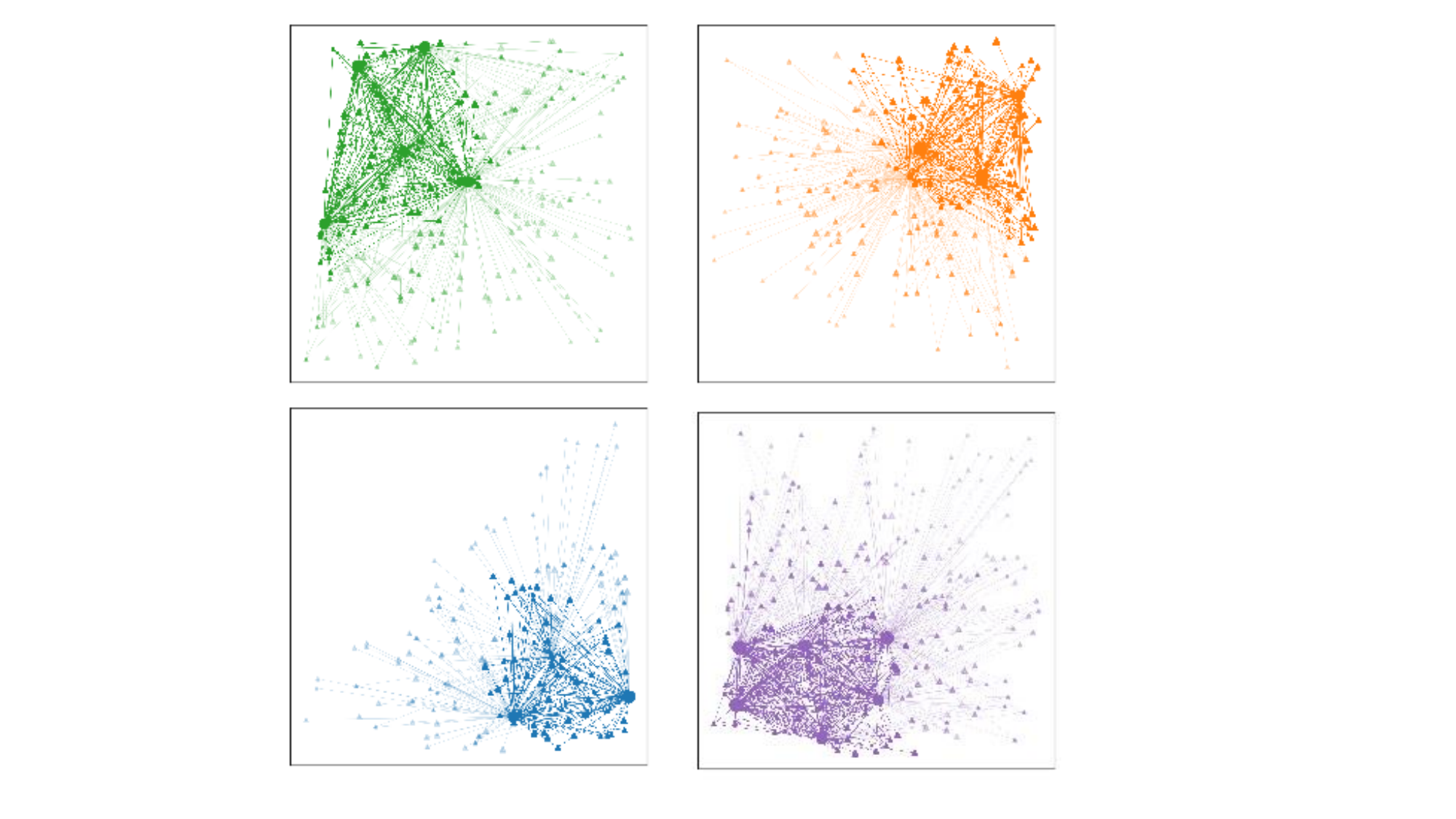}%
    }

    \caption{The initial kernel and the sequence of buckets produced in Phase~2 of \alg\ 
    for instance \texttt{i300\_1}.}
    \label{fig:KernelBucketsPKS}

\end{figure}
Figure~\ref{fig:InitKPKS} illustrates, for instance {\texttt{i300\_1}}, the composition of the initial kernel $K_y$ (the red circles) along with the set of customers grouped in regions. For readability, the assignment variables composing the initial kernel $K_x$ are not displayed. Figure~\ref{fig:BucketsPKS} depicts the composition of the four buckets determined by \alg. For each picture in Figure~\ref{fig:BucketsPKS}, lighter colors are used to highlight inter-regional allocations. For the same instance, Figure~\ref{fig:KernelBucketsKS2014} displays the initial kernel and the sequence of six buckets determined by \ks. It is evident that, by leveraging the regions derived in Phase 1, \alg\ produces significantly smaller buckets, excluding more $y$ and $x$ variables from the subsequent optimization than \ks. The results of the experimental analysis reported below highlight that these exclusions lead to an improved algorithmic efficiency, without deteriorating the quality of the solutions found. On the contrary, for numerous instances \alg\ significantly outperforms \ks\ in terms of quality of the solution produced.

Once the initialization step is concluded, Phase 2 of \alg\ continues along the general lines of the standard KS framework. This step is henceforth referred to as the \textit{optimization step}. Firstly, $BIP(K)$ is solved to obtain a first integer feasible solution, which is used to set an initial upper bound $z^H$. Subsequently, \alg\ seeks to improve this solution by solving a sequence of restricted BIPs. In each iteration, a modified restricted problem $BIP(K \cup B^r)$ is solved. An objective function cutoff value equal to the upper bound $z^H$ is always set. A second constraint is possibly added (see below for further details) to ensure that at least one binary variable from bucket $B^r_y$ is selected. By construction, if such a modified $BIP(K \cup B^r)$ is feasible, it improves upon the incumbent solution. In this case, kernel $K$ is modified by adding all new promising variables, and removing any that are no longer promising. After solving a number of restricted BIPs equal to the number of buckets generated, the search for an improving solution, and so the \alg, ends. Section~\ref{sec:Phase2} provides a detailed description of Phase 2.
 
\subsubsection{Phase 1: Pattern recognition.}\label{sec:Phase1}

In Phase 1, we gather information on the importance of and relationships between decision variables by solving a sequence of LP relaxations and extracting information via pattern recognition (see Algorithm~\ref{alg:PattRecogn}). To make the paper self-contained, we hereafter briefly describe spectral biclustering, the pattern recognition technique we applied, and highlight the foremost differences compared to the method applied in~\cite{Bakker2024}. The interested reader is referred to the latter article for further details.

Recall that $LP(\mathcal{Y} \cup \mathcal{X})$ denotes the LP relaxation of the \model, and let $s_1$ denote its optimal solution. Phase 1 begins solving $LP(\mathcal{Y} \cup \mathcal{X})$, followed by the solution of a sequence of restricted LP relaxations derived from $s_1$. 
As previously mentioned, inequalities \eqref{SSCFL-VI} are redundant, but have a strong impact on both the tightness of the LP relaxation (positive impact) and the computing times (negative for large instances). 
Given the potential negative effects on computing times, Phase 1 begins solving the LP relaxation of the \model\ without constraints \eqref{SSCFL-VI}. Then, based on the solution obtained and some pivotal characteristics of the instance, \alg\ decides either to enforce that a sufficient number of $y$ variables take a nonzero value or to iteratively add a subset of constraints \eqref{SSCFL-VI} to tighten the LP relaxation. The number of these iterations is data-driven and adaptive upon the following characteristics of the instance. Notice that a rough estimate of the number of locations open in an optimal solution to the \prob\ is $I^*=\vert I \vert / \rho$, where $\rho$ denotes the ratio between the total capacity and the total demand --i.e., $\rho = \frac{\sum_{i \in I} q_i}{\sum_{j \in J} d_j}$. Intuitively, the larger $\rho$, the larger the total capacity relative to the demand, and the smaller the number of locations required to serve all customers. Let $s^\prime$ be the optimal solution to the LP relaxation found in a given iteration, and denote by $I^\prime$ the number of $y$ variables taking a nonzero value in $s^\prime$. If $I^\prime$ is greater than $(1+\epsilon) \cdot I^*$, for a small $\epsilon>0$, for each variable $x_{ij}$ taking a nonzero value in $s^\prime$ we add the corresponding inequality \eqref{SSCFL-VI}. This process stops as soon as solution $s^\prime$ is such that satisfies $I^\prime \leq (1+\epsilon) \cdot I^*$ or a maximum number $\#VI$ of iterations has been carried out. Solution $s^\prime$ produced in the last iteration is used as the optimal solution to the LP relaxation $s_1$. On the other hand, if solving the LP relaxation of the \model\ without constraints \eqref{SSCFL-VI} the value $I^\prime$ is smaller than $I^*$, this often suggests that more locations, compared to the number of $y$ variables taking a nonzero value in solution $s^\prime$, are required to serve all customers in a feasible solution to the \prob. Hence, to enforce that more $y$ variables are selected in the LP relaxation, we add the constraint $\sum\nolimits_{i \in I} y_i \geq I^*$ and solve the new LP relaxation to obtain $s_1$. In both cases, solution $s_1$ is used to initialize set $\mathcal{S}$.  

Next, a sequence of $N$ restricted LP relaxations $LP(U_y \cup U_x)$ is solved to produce the remaining solutions. The process of generating this sequence is guided by the $y$ variables taking a nonzero value in solution $s_1$, as in each of the $N$ solutions $\alpha$ of these variables are fixed to 0 forcing others to be selected instead. More precisely, the number $\alpha$ of $y$ variables removed in each of the restricted problems is equal to $\lceil I^1 / \min \{\lceil \rho\rceil,N\} \rceil$, where $\lceil \cdot \rceil$ denotes the round up to the nearest integer of a given value. In other words, $\alpha$ is inversely related to ratio $\rho$ when the latter is smaller than $N$. The rationale is that, in instances with a small $\rho$, several locations will be open in an optimal solution to the \prob\ and, hence, more candidates must be considered in the restricted BIPs solved in Phase 2. In fact, having a large value of $\alpha$ implies a stronger perturbation of solution $s_1$ and enables the gathering of information on more alternative locations. Evidently, $\alpha$ cannot be as large as to render the resulting restricted LP relaxation infeasible. To avoid this situation, \alg\ checks whether the remaining $|I|-\alpha$ locations are expected to be sufficient to satisfy the total demand, and if this is not the case, $\alpha$ is reduced accordingly.

\begin{algorithm}
{\setlength{\baselineskip}{0.75\baselineskip} 
  \footnotesize 
\begin{algorithmic}[1]
\caption{\alg\ Phase 1: Pattern recognition}
\label{alg:PattRecogn}
\State \Comment{\textbf{Produce the first solution $s_1 \in \mathcal{S}$.}}
\State Let $LP(\mathcal{Y} \cup \mathcal{X})$ be the LP relaxation of the \model\ without inequalities \eqref{SSCFL-VI}.
\State Let $I^* = \vert I \vert / \rho$, where $\rho = \sum_{i \in I} q_i / \sum_{j \in J}d_j $.
\State Solve $LP(\mathcal{Y} \cup \mathcal{X})$. Let $s^\prime = \{ \bm{y}^\prime, \bm{x}^\prime \}$ be the optimal solution found, and $I^\prime = \vert \{i \in I \vert y_i'>0\} \vert$.
\If{$I^\prime < I^*$} \Comment{{Ensure enough $y$ variables take a nonzero value.}}
    \State Add $\sum_{i \in I} y_i \geq I^*$ to $LP(\mathcal{Y} \cup \mathcal{X})$ and solve the latter.
    \State Let $s_1$ be the optimal solution found.
\Else  \Comment{{Lazily add \eqref{SSCFL-VI}.}}
    \State $\ell \gets 1$. 
    \While{$I^\prime > (1+\epsilon) \cdot I^*$  \textbf{and} $\ell \leq  \#VI$}
        \State For each pair $(i^\prime,j^\prime)$ such that $x_{i^\prime j^\prime} > 0$ in $s^\prime$, add to $LP(\mathcal{Y} \cup \mathcal{X})$ inequality $x_{i^\prime j^\prime} \leq y_{i^\prime}$.
        \State Solve $LP(\mathcal{Y} \cup \mathcal{X})$. Let $s^\prime$ be the optimal solution found, and update $I^\prime$ accordingly.
        \State $\ell \gets \ell+1$.
    \EndWhile
    \State Let $s_1$ be the optimal solution found in the last iteration.
\EndIf
\State Initialize set $\mathcal{S}$ with solution $s_1 = \{ \bm{y}^1, \bm{x}^1 \}$, and let $I^1 = \vert \{i \in I \vert y^1_i > 0\}\vert$.
    \State \Comment{\textbf{Produce solutions $s_2, \dots, S_{N+1} \in \mathcal{S}$.}}
    \State Let $\alpha = \lceil I^1 / \min \{\lceil \rho\rceil,N\} \rceil$ be the number of $y$ variables to remove. \Comment{Determine perturbation level.}
    
    \State If removing $\alpha$ variables is expected to make the problem infeasible, reduce $\alpha$ accordingly.
    \For{$n \gets 1$ \textbf{to} $N$}
    \State Randomly select $\alpha$ of the $y$ variables taking nonzero value in $s_1$. Let $U_{y^n}$ be the set of $y$ variables not selected.
    \State Solve $LP(U_{y^n} \cup U_x)$ to obtain $s_{n+1}$ and add the latter to $\mathcal{S}$.
    \EndFor
\State \Comment{\textbf{Derive the set of regions $\mathcal{R}$ from $\mathcal{S}$.}}
\State Generate feature matrix $A$, with entries according to \eqref{eq:AEntries}, and determine $I_0$.
\State $r \gets 1$, $\mathcal{R}\gets\{(I \backslash I_0,J)\}$
\While {$\ell^{inter}(\mathcal{R},A)<\theta$}
\State  $r \gets r+1$
\State  Derive $r$ regions via spectral biclustering, such that $\mathcal{R}\gets\{(I_1,J_1),\dots,(I_{r},J_{r})\}$.
\EndWhile
\end{algorithmic}
}
\end{algorithm}

The solutions in $\mathcal{S}$ serve as an input to the biclustering algorithm that groups locations and customers into a set of regions $\mathcal{R}$. 
Biclustering is a clustering technique that partitions simultaneously the rows and columns of a feature matrix \citep[see, among others,][]{busygin2008biclustering}. In this paper, we use a spectral biclustering algorithm, which is a matrix factorization-based technique that uses the information from a feature matrix to derive a bipartite undirected graph in which the edge weights are the matrix entries. Then, the technique solves a partitioning problem via spectral graph partitioning \citep[see, for further details,][]{dhillon2001co}. As many conventional clustering techniques, it requires as an input the number of biclusters to produce and then strives for a high intra-bicluster similarity and a high inter-bicluster dissimilarity. 

To derive a feature matrix $A$ from the set of solutions $\mathcal{S}$, we aggregate the allocation matrices of the individual solutions as follows:  
\begin{equation}\label{eq:AEntries}
a_{ij} = \sum_{s=1}^{\vert \mathcal{S}\vert} \left\lceil x^{s}_{ij} \right\rceil \quad (i \in I; j \in J),    
\end{equation}
where $\left\lceil x^{s}_{ij} \right\rceil$ denotes the value taken by variable $x_{ij}$ in solution $s \in \mathcal{S}$ (i.e., $x^{s}_{ij}$) rounded up to the nearest integer. In other words, each entry $a_{ij}$ of $A$ equals the number of times customer $j$ is served, even with a small value, from location $i$ across all solutions in $\mathcal{S}$. This stands in contrast with \cite{Bakker2024}, where all entries of $A$ indicated whether or not an assignment was nonzero in any of the solutions in $\mathcal{S}$.
Let $I_0$ denote the set of locations that have never been selected in any solution in $\mathcal{S}$. Notice that each row in $A$ corresponding to a location in $I_0$ contains only zero entries, and thereby provides no information on the subset of customers these locations may serve. Thus, as previously mentioned, these locations are removed from the optimization. In terms of the biclustering algorithm, this implies that $A$ is reduced to a $\vert I \backslash I_0 \vert \times \vert J \vert$-dimensional matrix, and that the output of the biclustering is a set of regions $\mathcal{R}$ that does not include any location from $I_0$.

The purpose of applying biclustering and deriving regions is to later generate the initial kernel and the buckets so that they capture as much of the relevant interactions between the decision variables as possible. Hence, we seek to derive regions in such a way that there is a strong interaction between the locations and customers within the region, and little interaction between locations and customers of different regions. To this aim, we capture inter-cluster interaction using the fraction of out-of-region allocations across solutions in $\mathcal{S}$ defined in \cite{Bakker2024} as the single internal validation criterion. The number of out-of-region allocations in proportion to the total number of allocations across solutions in $\mathcal{S}$ is determined as:
$$\ell^{inter}(\mathcal{R},A) = \frac{\sum\limits_{r=1}^{\vert \mathcal{R}\vert} \sum\limits_{i \in I_r} \sum\limits_{j \in J \backslash J_r} a_{ij}}{\sum\limits_{i \in I} \sum\limits_{j \in J} a_{ij}}.$$
As noted in \cite{Bakker2024}, regions implied in an instance of a location problem are an abstract concept, and there is no true number of regions underlying a problem. Instead, we look for the maximum number of regions that can be identified using spectral biclustering without exceeding an upper threshold $\theta \in [0, 1]$ on $\ell^{inter}(\mathcal{R},A)$. Therefore, starting from one cluster comprising all locations in $I \backslash I_0$ and customers in $J$, \alg\ repeatedly increases the target number of clusters, denoted as $r$, until the threshold is exceeded.

\subsubsection{Phase 2: Kernel Search.}\label{sec:Phase2}

\begin{algorithm}
{\setlength{\baselineskip}{0.75\baselineskip} 
  \footnotesize 
\begin{algorithmic}[1]
\caption{\alg\ Phase 2: Initialization step}
\label{alg:PKSInitalization}
\State \Comment{\textbf{Compute a score for each location.}}
\State Assign weight $\lambda_s \geq 0$ to each solution $s \in \mathcal{S}$.
\State For each location $i \in I$, compute score $w_i$ as in \eqref{eq:scoresW}.
\State \Comment{\textbf{Determine initial kernel $K_y$ and buckets $B^r_y$.}}
\For{ $r \gets 1$ \textbf{to} $|\mathcal{R}|$}
\State Sort the $y$ variables associated with the locations in region $r$ in non-increasing order of their score $w_i$.
\State Include in the initial kernel $K_y$ the first $m_r$ variables.  \Comment{Assign the most promising to $K_y$.}
\State Assign each remaining $y$ variable, with nonzero score $w_i$, to bucket $B^r_y$.
\EndFor
\State Sort the sequence of buckets $\{B^r_y\}_{r=1,\dots,|\mathcal{R}|}$ in non-increasing order of their cardinality.
\State \Comment{\textbf{Determine initial kernel $K_x$.}}
\For{each $y_i \in K_y$}
\State Include in $K_x$ all $x_{ij}$ variables such that:
\State $\quad$ (1) $i$ and $j$ are both in the same region and $\hat{c}(x_{ij}) \leq \hat{\gamma}$. \Comment{Add intra-regional allocations.}
\State $\quad$ (2) $i$ and $j$ are in different regions, $\hat{c}(x_{ij}) \leq \hat{\gamma}$, and $c_{ij} \leq \gamma_j$. \Comment{Add gray zone allocations.}
\EndFor 
\State \Comment{\textbf{Determine buckets $B^r_x$.}}
\For{ $r \gets 1$ \textbf{to} $|\mathcal{R}|$}
\State For each $y_i \in B^r_y$, include in $B^r_x$ all $x_{ij}$ variables representing intra-regional and gray zone allocations.
\EndFor
\end{algorithmic}
}
\end{algorithm}
Phase 2 of \alg\ has two main steps. First, the initialization step outlined in Algorithm~\ref{alg:PKSInitalization} is carried out. It begins assigning a non-negative weight $\lambda_s$ to each solution $s \in \mathcal{S}$. The basic rationale of these weights is to prioritize those LP relaxations that are the most likely to provide accurate indications on which variables are promising. Recall that in Phase 1 \alg\ solves, first, the LP relaxation of the original problem, and then a sequence of restricted LP relaxations. Hence, one might expect that an optimal solution to the former can provide better indications than the latter restricted LP relaxations. Next, for each location $i \in I$ a non-negative score $w_i$ is computed as a linear weighted combination of the amount of demand served across all solutions in $\mathcal{S}$. More formally, for each location $i \in I$ the score is computed as follows:
\begin{equation}\label{eq:scoresW}
    w_i = \sum_{s = 1}^{N+1}\lambda_{s}\sum_{j\in J} d_j x_{ij}^{s} \quad (i \in I).
\end{equation}
Note that the use of weights $\bm{\lambda}$ and scores $\bm{w}$, as defined above, generalizes \ks, which, in this regard, is a special case obtained by setting $\bm{\lambda} = (1, 0, \dots, 0)$. Further, notice that, by construction, $w_i$ is equal to zero for each location $i$ in set $I_0$. In contrast with the standard KS framework, in \alg\ the composition of the initial kernel, as well of the sequence of buckets, is guided by the following key factors: the set of optimal LP relaxation solutions $\mathcal{S}$, the scores $\bm{w}$, and the set of regions $\mathcal{R}$. More precisely, for each region $r \in \mathcal{R}$, \alg\ determines a ranking $L_r$ of the corresponding locations by sorting them in non-increasing order of their score $w_i$. The initial composition of kernel $K_y$ is determined as follows. Let $m_r$ be the number of $y$ variables associated with the locations in region $r$ that take a nonzero value in solution $s_1$. Then, kernel $K_y$ is initialized by taking, for each region $r$, the $y$ variables associated with the first --and hence, most promising-- $m_r$ locations in ranking $L_r$. Next, an initial sequence of buckets is determined creating one bucket $B^r_y$ for each region $r \in \mathcal{R}$. Each of these buckets is taken to be the set comprising each $y_i$ variable such that: (i) $y_i$ has not been added already to $K_y$, (ii) location $i$ belongs to region $r$, and (iii) score $w_i$ is nonzero. Subsequently, these buckets are sorted in non-increasing order of their cardinality producing the sequence $\{B^r_y\}_{r=1,\dots,|\mathcal{R}|}$. Note that some buckets can possibly be empty and, after the sorting, are in the last positions of the sequence. To simplify the exposition, we hereafter assume that all buckets are non-empty. If, on the contrary, some buckets are empty, each procedure described in the following that explores through the sequence of buckets stops after the last non-empty one. 

Kernel $K_x$ is initialized by adding an appropriately chosen subset of the $x$ variables in two consecutive steps. Let $\hat{c}(x_{ij})$ be the reduced cost coefficient of variable $x_{ij}$ in solution $s_1$. Additionally, let $\hat{\gamma}$ be a threshold on the reduced cost coefficients used to select the $x$ variables. Then, for each location $i$ such that $y_i \in K_y$, in the first step \alg\ adds to the initial kernel $K_x$ all $x_{ij}$ variables such that: (i) customer $j$ belongs to the same region as location $i$, and (ii) $\hat{c}(x_{ij}) \leq \hat{\gamma}$. Note that these pairs $(i,j)$ represent intra-regional allocations, and one of them can be expected to be selected in an optimal solution to the original problem. Akin to \ks, we compute the threshold $\hat{\gamma}$ as the median of the reduced cost coefficients of the $x$ variables in solution $s_1$. After preliminary experiments, we realized that adding only the intra-regional allocations is too restrictive, as sometimes customers that are located near the boundary of a given region are, in an optimal solution to the original problem, allocated to a location belonging to a different, usually adjacent, region. These pairs $(i,j)$ represent inter-regional allocations. We call \textit{gray zones} all the areas comprising the customers that might be served from locations belonging to a different region. Since the idea of gray zones is to capture optimal assignments to customers across regional borders yet close to the border, we identify the associated $x_{ij}$ variables as follows. Let $\gamma_j$ be a threshold on the value of parameter $c_{ij}$. Then, for each location $i$ such that $y_i \in K_y$, in the second step \alg\ adds to the initial kernel $K_x$ all $x_{ij}$ variables such that: (i) customer $j$ belongs to a different region than location $i$, (ii) $\hat{c}(x_{ij}) \leq \hat{\gamma}$, and (iii) $c_{ij} \leq \gamma_j$. After preliminary experiments, for each customer $j \in J$ the threshold $\gamma_j$ is computed as the median of the values of $c_{ij}$.  Subsequently, each bucket $B^r_x$ is determined by selecting, for each $y_i \in B^r_y$, all $x_{ij}$ variables that fulfill the intra-regional and gray zone allocations criteria described previously for the initial kernel composition. The output of the initialization step is the initial kernel $K = K_y \cup K_x$ and the sequence of buckets $\{B^r = B^r_y \cup B^r_x\}_{r=1,\dots,\vert \mathcal{R} \vert}$. Note that some variables might not have been included neither in the initial kernel, nor in any bucket. Hence, they are excluded from the optimization by permanently fixing their value to zero in any restricted BIP. Indeed, recall that \alg\ fixes to zero each $y_i$ such that $i \in I_0$. Similarly, each $x_{ij}$ variable that does not meet neither the intra-regional nor gray zone allocations criteria is fixed to zero as well. 
\begin{algorithm}
{\setlength{\baselineskip}{0.75\baselineskip} 
  \footnotesize 
\begin{algorithmic}[1]
\caption{\alg\ Phase 2: Optimization step}
\label{alg:PKSOptimization}
\State \Comment{\textbf{Produce an initial feasible solution.}}
\State Solve $BIP(K)$.
\While{$BIP(K)$ is infeasible }
\State Add to $K$ the following bucket $B^r$.
\State Solve $BIP(K)$.
\EndWhile

\If{$BIP(K)$ is feasible}
\State Initialize $z^H$ and $(\bm{y}^H, \bm{x}^H)$ with the best solution found.
\Else
\State \textbf{return} {\ttfamily{fail}} \Comment{Time limit is reached before an initial feasible solution was found.}
\EndIf

\State \Comment{\textbf{Iterate over the remaining buckets.}}
\While{not all buckets have been analyzed}
\State Solve modified $BIP(K \cup B^r)$.
\If{modified $BIP(K \cup B^r)$ is feasible}\Comment{Learn and adjust.}
\State Add to $K_y$ each $y_i$ variable from $B^r_y$ taking value 1, and add to $K_x$ all the associated $x_{ij}$ from $B^r_x$.
\State Remove from $K_y$ each $y_i$ not selected in any of the $p$ previous iterations, and remove from $K_x$ all the associated $x_{ij}$.
\State Update $z^H$ and $(\bm{y}^H, \bm{x}^H)$.
\EndIf
\EndWhile
\end{algorithmic}
}
\end{algorithm}

\alg\ carries on executing the optimization step sketched in Algorithm~\ref{alg:PKSOptimization}. Model $BIP(K)$ is solved first. Assuming that the original problem is feasible, if $BIP(K)$ is infeasible,  the initial kernel is likely too small. In these cases, kernel $K$ is iteratively adjusted by adding the variables from the following bucket in the sequence until problem $BIP(K)$ is feasible, or the total time limit is reached. In the latter case, \alg\ fails to produce a feasible integer solution to the \prob. In the former case, \alg\ initializes the upper bound $z^H$ and the best solution found $(\bm{y}^H, \bm{x}^H)$. Next, a sequence of restricted BIPs, one for each of the remaining buckets, is solved. In each iteration, a restricted $BIP(K \cup B^r)$ model is solved after the introduction of (possibly) both the following constraints. First, the upper bound $z^H$ is always introduced to set an objective function cutoff value, so that any solution with a value worse than or equal to $z^H$ is discarded. The second constraint of the form $\sum\nolimits_{i | y_i \in B^r_y} y_i \geq 1$ ensures that at least one location of the new bucket $B^r_y$ is selected. Here, it is worth pointing out that while solving large-scale instances restricted $BIP(K \cup B^r)$ models are often not solved to proven optimality within the time available. To handle these situations, we proceed as follows. If a restricted BIP, say $BIP(K \cup B^{r-1})$, is solved to optimality, to improve upon the incumbent solution any feasible solution to the following restricted BIP --i.e., $BIP(K \cup B^r)$-- must include at least one $y$ variable from the bucket associated. Only in this case we add the second constraint mentioned above.


By construction, if this modified $BIP(K \cup B^r)$ is feasible, it improves upon the incumbent solution. The composition of kernel $K$ is modified accordingly. Each variable from bucket $B^r_y$ that takes value 1 in the improving solution is added to kernel $K_y$. Accordingly, kernel $K_x$ is adjusted by including all the associated assignment variables from bucket $B^r_x$. To control the size of kernel $K$, each $y$ variable that has not been selected in any of the previous $p$ iterations is removed from the kernel, along with the associated $x$ variables. Finally, the upper bound $z^H$ and the incumbent solution are updated.  

\alg\ ends when all restricted BIPs have been solved. To control the computing time spent for solving large-scale instances , a limit to the latter is also imposed. In detail, at the beginning of the optimization step, the time remaining after carrying out Phase 1 and the initialization step is equally divided for the solution of the $1 + \vert \mathcal{R} \vert$ restricted BIPs. If, after solving any of the latter BIPs, some of the time initially allotted has not been used, the remainder is proportionally added to the time assigned for the solution of each of the following restricted BIPs. 

\section{Experimental analysis}\label{sec:Experiments}

This section is devoted to the presentation and discussion of the computational experiments. It is organized as follows. Section~\ref{sec:Testing} describes the hardware and software used, along with the data sets, and the control parameters of \alg. In Section~\ref{sec:AgainstLiterature}, we validate, on a set of benchmark instances, the performance of \alg\  against \cpx, a re-implementation of \ks, and \hils. To ensure comparability with \alg, we re-implemented \ks, hereafter simply termed \ks, so that both methods are coded in the same programming language and use the same version of \cpx.  Finally, the results of the computational experiments on new very large-scale instances are reported in Section~\ref{sec:NewInstances}, where the performance of \alg\ is benchmarked against \cpx\ and \ks.

\subsection{Testing environment and implementation details}\label{sec:Testing}

The computational experiments were conducted on the high-performance computing cluster bwUniCluster 2.0\footnote{\url{https://www.scc.kit.edu/dienste/bwUniCluster_2.0.php} (last accessed: 21st of March 2025).}. Each instance of the computing cluster was allocated 64 GB of RAM, along with 20 CPU cores, each of the latter with a processor frequency of 2.1 GHz. According to extensive preliminary experiments, these settings ensure an average performance comparable to that of a current high-end workstation. 

The \model, \ks, and \alg\ have all been implemented in Python 3.8. Mathematical models were formulated using IBM’s DOcplex library (version 2.23.222) and solved with \cpx\ (version 22.11). After preliminary experiments, other than a limit on the computing time (details are below), \cpx\ parameters were left at their default values.

In the computational experiments, we used two data sets. The first data set is composed of benchmark instances for the \prob\ taken from the literature, whereas the second comprises a subset of very large-scale benchmark instances originally proposed for the \mscflp. Table~\ref{tab:inst} reports a summary of the main characteristics of the instances composing the test sets. Notice that several authors acknowledge that the ratio $\rho$ between the total capacity and the total demand has a pivotal impact on the performance of solution algorithms for capacitated location problems \citep[see, among others, ][]{Klose2007, Avella2008}. Thus, and also given its impact on how \alg\ functions, in the rightmost column of Table~\ref{tab:inst} we report the range of values taken by $\rho$.

\begin{table}
\centering
\caption{A summary of the instances in the test sets.}
\label{tab:inst}
\begin{tabular}{@{}l@{\quad}l@{\quad}c@{}@{\quad}r@{\quad}r@{\quad}c}
\hline
 & \textbf{Data Set}              & \textbf{\# Inst.} & $\bm{|I|}$  & $\bm{|J|}$  & $\bm{\rho}$         \\ \hline 
\parbox[t]{2mm}{\multirow{6}{*}{\rotatebox[origin=c]{90}{\textbf{Benchmark}}}} & {\ttfamily{OR4}}    & 12 & 100 & 1,000 & $9.7 \leq \rho \leq 27.5$       \\ 
\cline{2-6}
& {\ttfamily{TB1}} & 20  & 300 & 300 & \\
& {\ttfamily{TB2}} & 20  & 300 & 1,500 & \\
  & {\ttfamily{TB3}} & 20  & 500 & 500 & $\rho \in \{5, 10, 15, 20\}$\\
& {\ttfamily{TB4}} & 20  & 700 & 700 & \\
& {\ttfamily{TB5}} & 20  & 1,000 & 1,000 & \\
\hline
\hline
\parbox[t]{2mm}{\multirow{5}{*}{\rotatebox[origin=c]{90}{\textbf{New}}}} &  {\ttfamily{TB-A1/B1/C1}} & 30/25/30  & 800 & 4,400  & \\
& {\ttfamily{TB-A2/B2/C2}} & 25/29/30  & 1,000 & 1,000 & \\
& {\ttfamily{TB-A3/B3/C3}} & 30/30/30  & 1,000 & 4,000 & $\rho \in \{1.1, 1.5, 2, 3, 5, 10\}$\\
& {\ttfamily{TB-A4/B4/C4}} & 30/30/30  & 1,200 & 3,000 &  \\
& {\ttfamily{TB-A5/B5/C5}} & 25/25/26  & 2,000 & 2,000 &  \\ \hline
\end{tabular}
\end{table}

The first test set is composed of 112 benchmark instances for the \prob\ used as a test bed by other authors, including \citet{guastaroba2014heuristic}, \citet{caserta2020general}, and \citet{de2024hybrid}. Recall that the latter paper introduced \hils, which, to our knowledge, is the current state-of-the-art algorithm for solving the \prob. These instances are grouped into two subsets. The first, henceforth denoted as {\ttfamily{OR4}}, is composed of the largest instances currently belonging to the OR-Library, and consists of 100 locations and 1,000 customers\footnote{\url{https://www.brunel.ac.uk/~mastjjb/jeb/orlib/capinfo.html} (last accessed: 21st of March 2025).}. The second subset is composed of the 100 instances tested in \citet{Avella2009} for the \mscflp\ that turned out to be feasible for the \prob\ according to the results provided in \citet{guastaroba2014heuristic}\footnote{\url{https://or-brescia.unibs.it/instances/instances_sscflp} (last accessed: 21st of March 2025).}. These instances, hereafter denoted by {\ttfamily{TB1}} through {\ttfamily{TB5}}, range from 300 to 1,000 locations, and from 300 to 1,000 customers. We used this set of benchmark instances to validate the performance of \alg\ against \cpx, \ks, and \hils. Each of the benchmark instances is solved with the latter two approaches and \alg, and each of the three methods was run with a limit on the total computing time of 3,600 seconds. The results of the validation on the benchmark instances are discussed in Section~\ref{sec:AgainstLiterature}.

The second test set comprises 425 very large-scale instances. These are a subset of the 445 instances originally proposed for the \mscflp\ in \citet{Avella2008} and \citet{guastaroba2012kernel}, and not all are necessarily feasible for the \prob. To identify the feasible instances, we set a time limit of 7,200 seconds and solved with \cpx\ each of the original instances. Eventually, we discarded 20 instances where \cpx\ did not find any feasible solution within the time limit imposed. The 425 remaining instances are clustered according to the ratio between the average opening and assignment costs into three main data sets denoted by {\ttfamily{TB-A}}, {\ttfamily{TB-B}}\footnote{\url{http://wpage.unina.it/sforza/test/} (last accessed: 21st of March 2025).}, and {\ttfamily{TB-C}}\footnote{\url{https://or-brescia.unibs.it/instances/instances_clfp} (last accessed: 21st of March 2025).}. Each of the latter three data sets is further classified according to the size of the instances (see Table~\ref{tab:inst}). These instances range from 800 to 2,000 locations, and from 1,000 to 4,400 customers. As these are novel instances for the \prob, we validate the performance of \alg\ against \cpx\ and \ks. Each of the large-scale instances was solved with each of these three approaches imposing a limit on the computing time equal to 7,200 seconds. 

Table~\ref{tab:InstProblems} summarizes the issues that each method encountered while solving the new large-scale instances. In detail, column \textit{\#Fails} indicates the number of instances a method did not find a feasible integer solution to the \prob. For \cpx, we provide also column \textit{\#OoM} that gives the number of times the solver terminated its search before reaching the time limit due to an out of memory error. Looking at these values, it is evident that \cpx\ ran out of memory for roughly half of the instances. In these cases, we recorded the best solutions and best lower bounds found by \cpx\ before crashing. However, comparing the computing times of \cpx\ with the other methods is misleading, since the error usually occurred well before the time limit. For this reason, we decided not to include the computing times of \cpx\ in the analysis provided in Section~\ref{sec:NewInstances}. Notice that \alg\ is the approach that encountered the smallest number of problems. This emphasizes the reliability and robustness of the approach proposed in the present paper and, after observing the number of failures of the standard \ks, the beneficial impact of the methodological contributions we introduced. Indeed, \ks\ did not find a feasible integer solution for 86 of the 425 large-scale instances. To ensure consistency and fair comparability of the experimental results, the performance evaluation provided in Section~\ref{sec:NewInstances} covers only the 337 instances where each method found a feasible integer solution (see column \textit{\#All Feas.}). The detailed computational results for all the new large-scale instances are available in the online repository.


\begin{table}[htbp]
\centering
\caption{New large-scale instances: Summary of issues encountered by each approach.}
\label{tab:InstProblems}
\begin{tabular}{@{\extracolsep{4pt}}lc rr r r r}
\hline 
&  &    \multicolumn{2}{l}{\textbf{\cpx}}  & \multicolumn{1}{l}{\textbf{\ks}} &   \multicolumn{1}{l}{\textbf{\alg}}  \\ 
             \cline{3-4} \cline{5-5} \cline{6-6}
   \textbf{Data }              &  &      &  &      &      &  \multicolumn{1}{c}{\textbf{\#All}} \\
  \textbf{Set}              & \textbf{\#Inst.} & \multicolumn{1}{c}{\textbf{\#Fails}}     & \multicolumn{1}{c}{\textbf{\#OoM}} &  \multicolumn{1}{c}{\textbf{\#Fails}}    &  \multicolumn{1}{c}{\textbf{\#Fails}}     &  \multicolumn{1}{c}{\textbf{Feas.}} \\
  \hline 
{\ttfamily{TB-A}} & 140 & 0 & 105 & 36  & 0  & 104  \\ 
{\ttfamily{TB-B}} & 139 & 0 & 81 & 33  & 3  & 104  \\
{\ttfamily{TB-C}} & 146 & 0 & 20 & 17  & 1  & 129   \\
\cline{1-2} \cline{3-4} \cline{5-5} \cline{6-6} \cline{7-7} 
\textbf{Tot.} & \textbf{425}  & \textbf{0} & \textbf{206} & \textbf{86}  & \textbf{4}  & \textbf{337}   \\ \hline
\end{tabular}
\end{table}

Finally, after extensive preliminary experiments, we chose the configuration of the parameters controlling \alg\ reported in Table~\ref{tab:param}.

\begin{table}
\centering
\caption{\alg: Best parameters configuration}\label{tab:param}
\begin{tabular}{@{}l@{\quad}l@{}@{\quad}r}
\hline 
\textbf{Parameter}              & \textbf{Description} & \textbf{Value}            \\ \hline
$N + 1$    & Phase 1: Number of LP relaxations solved & 10 + 1       \\ 
$\#VI$ & Phase 1: Maximum number of iterations of adding inequalities \eqref{SSCFL-VI}  & 5 \\
$\theta$ &  Phase 1: Threshold on inter-regional validity & 0.05 \\ 
\hline
$\bm{\lambda}$   & Phase 2: Weights assigned to solutions in $\mathcal{S}$ & $(1, N^{-1}, \dots, N^{-1})$  \\
$p$ & Phase 2: Number of iterations before removing variables from kernel $K$ & 2 \\ \hline
\end{tabular}
\end{table}

\subsection{Performance evaluation on the benchmark instances} \label{sec:AgainstLiterature}

In this section, we validate the performance of \alg\ on the benchmark instances. The quality of a solution is compared against the optimal solution (if available) or the best upper bound determined as follows. As mentioned above, each of the benchmark instances was solved, among others, in \citet{guastaroba2014heuristic}, \citet{caserta2020general}, and \citet{de2024hybrid}, where the authors provided detailed results. In addition, we solved each instance with \cpx, \ks, and \alg. For each instance, we determined the best solution value, denoted by $z^{UB}$, comparing the results reported in the three aforementioned papers with those we obtained. In the following, for each method we provide statistic \textit{\#Best} that gives the number of instances on which the method found the best solution. Statistic \textit{$z^{UB}$ Gap (\%)} measures the deterioration with respect to $z^{UB}$. This statistic is computed as $100 \times \frac{z^H - z^{UB}}{z^{UB}}$, where $z^{H}$ is the value of the best solution produced by the method under consideration. Finally, we show for each method the computing times, measured in seconds and denoted by \textit{Time (sec.)}. An important remark is related to the results we considered for \hils. For each benchmark instance, the authors of \citet{de2024hybrid} ran five times \hils, and provided the best solution value found across the five runs, along with the average solution value and the average time. To compute the statistics above for \hils\ we used the best solution value indicated in \citet{de2024hybrid} instead of the average value, thereby being in favor of their approach. Regarding the computing times, to ensure comparability, we multiplied the average values indicated in \citet{de2024hybrid} by five--i.e., the number of runs performed to produce the best solution. The reader should be aware that the comparison on computing times has to be done with care, since \hils\ was run on a different hardware and coded in a different programming language, namely C++ that is well-known to be generally more efficient than Python.

 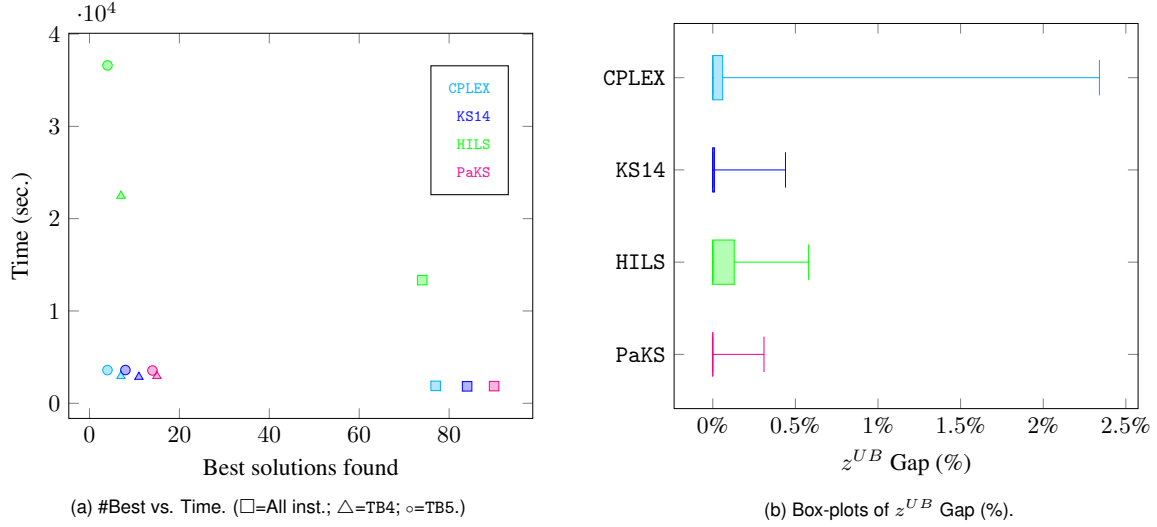
\begin{figure}
     \begin{minipage}[t]{0.5\textwidth}
\centering
\begin{tikzpicture}[scale=0.9]
\begin{axis}[
  ylabel={Time (sec.)},    
  xlabel={Best solutions found},
    scatter/classes={
        a1={mark=square*,cyan,fill opacity=0.3},
        a2={mark=square*,blue,fill opacity=0.3},
        a3={mark=square*,green,fill opacity=0.3},
        a4={mark=square*,magenta,fill opacity=0.3}, 
        b1={mark=triangle*,cyan,fill opacity=0.3},
        b2={mark=triangle*,blue,fill opacity=0.3},
        b3={mark=triangle*,green,fill opacity=0.3},
        b4={mark=triangle*,magenta,fill opacity=0.3},
        c1={mark=*,cyan,fill opacity=0.3},
        c2={mark=*,blue,fill opacity=0.3},
        c3={mark=*,green,fill opacity=0.3},
        c4={mark=*,magenta,fill opacity=0.3}%
    },
] 
\addplot[
        scatter, 
        only marks,
        scatter src=explicit symbolic,
    ]
    table[meta=label] {
        x       y       label   annotation
        77      1887    a1      \cpx
        84      1834    a2      \ks
        74      13340   a3      \hils
        90      1858    a4      \alg
    };
    \addplot[
        scatter, 
        only marks,
        scatter src=explicit symbolic,
    ]
    table[meta=label] {
        x       y       label   annotation
        7      2947    b1      \cpx
        11      2855    b2      \ks
        7      22467   b3      \hils
        15      2950    b4      \alg
        4      3595    c1      \cpx
        8      3602    c2      \ks
        4      36606   c3      \hils
        14      3550    c4      \alg
    };
\end{axis}
\matrix [draw, above left] at (6.5,3.3) {
  \node[cyan,fill opacity=0.9,font=\tiny] { \cpx}; \\
  \node[blue,fill opacity=0.9,font=\tiny] {\ks}; \\
  \node[green,fill opacity=0.9,font=\tiny] { \hils }; \\
    \node[magenta,fill opacity=0.9,font=\tiny] { \alg }; \\
};
     \end{tikzpicture}
     \subcaption{\#Best vs. Time. ($\square$=All inst.; $\triangle$={\ttfamily{TB4}}; $\medcirc$={\ttfamily{TB5}}.)}\label{fig:ScattBenchA}
\end{minipage}\hfill
\begin{minipage}[t]{0.5\textwidth}
   \begin{tikzpicture}[scale=0.9]
      \centering
  \begin{axis}
    [
    boxplot/draw direction=x,
    boxplot/box extend=0.48,
    ylabel={},
    xlabel={$z^{UB}$ Gap (\%)},
    cycle list={{magenta},{green},{blue},{cyan}},
    ytick={1,2,3,4},
    yticklabels={\alg, \hils, \ks, \cpx},
    xticklabel={\pgfmathparse{\tick}\pgfmathprintnumber{\pgfmathresult}\%} 
    ]
    \addplot+[
    fill,fill opacity=0.3,
    boxplot prepared={
      median=0.00,
      upper quartile=0.00,
      lower quartile=0.00,
      upper whisker=0.31, 
      lower whisker=0.00
    },
    ] coordinates {};
    \addplot+[
    fill,fill opacity=0.3,
    boxplot prepared={
      median=0.00,
      upper quartile=0.13,
      lower quartile=0.00,
      upper whisker=0.58, 
      lower whisker=0.00
    },
    ] coordinates {};
    \addplot+[
    fill,fill opacity=0.3,
    boxplot prepared={
      median=0.00,
      upper quartile=0.01,
      lower quartile=0.00,
      upper whisker=0.44, 
      lower whisker=0.00
    },
    ] coordinates {};
    \addplot+[
    fill,fill opacity=0.3,
    boxplot prepared={
      median=0.00,
      upper quartile=0.06,
      lower quartile=0.00,
      upper whisker=2.34, 
      lower whisker=0.00
    },
    ] coordinates {};
    \end{axis}
    \end{tikzpicture}
    \subcaption{Box-plots of $z^{UB}$ Gap (\%).}\label{fig:ScattBenchB}
\end{minipage}
\caption{Benchmark instances: A comparison of \cpx, \ks, \hils, and \alg.}\label{fig:ScattBench}
\end{figure}

To facilitate the interpretation of the results, Figure~\ref{fig:ScattBench} summarizes the computational results for the benchmark instances. The scatter plot in Figure~\ref{fig:ScattBenchA} shows, for each method, the value of \textit{\#Best} ($x$-axis) against the average computing time ($y$-axis), calculated across all instances considered. Each method is depicted with a different color, and squares denote the values computed across all benchmark instances. We also show the values for the instances in data sets {\ttfamily{TB4}} (triangles) and {\ttfamily{TB5}} (circles). To simplify the readability of the figure, we did not show the remaining data sets as the values of \textit{\#Best} across the four methods are very similar. The only notable difference is related to the average computing times of \hils\ compared to the other three approaches (see Table~\ref{tab:AgainstLiterature} for more details). The box-and-whisker plots in Figure~\ref{fig:ScattBenchB} show the distribution of \textit{$z^{UB}$ Gap (\%)}. The average value of each statistic broken down for each data set is detailed in Table~\ref{tab:AgainstLiterature}. 

The main insights we can gain from Figure~\ref{fig:ScattBench} and Table~\ref{tab:AgainstLiterature} are as follows:

\begin{itemize}
    \item[\textit{(a)}] \alg\ outperforms the other methods in terms of quality of the solution found:
    \begin{itemize}
        \item[\textit{(i)}] The number of best solutions found (\textit{\#Best}) by \alg\ is considerably larger than the other methods: in 90 instances out of 112 \alg\ produced the best solution, followed by \ks\ (84), \cpx\ (77), and \hils\ (74);
        \item[\textit{(ii)}] \alg\ achieved the smallest average \textit{$z^{UB}$ Gap (\%)} over the complete set of benchmark instances: 0.02\% against a slightly worse average gap of 0.03\% achieved by \ks, and followed by \hils\ (0.05\%), and \cpx\ (0.12\%);
        \item[\textit{(iii)}] The box-plots in Figure~\ref{fig:ScattBenchB} highlight that the range of values of \textit{$z^{UB}$ Gap (\%)} for the solutions produced by \alg\ is considerably smaller than that of all the other approaches, followed by \ks\ (the second-best method);
        \item[\textit{(iv)}] In terms of the average value of \textit{$z^{UB}$ Gap (\%)}, dissimilarities between the different approaches are noticeable only for the instances in data sets {\ttfamily{TB4}} and {\ttfamily{TB5}}. In both cases, \alg\ outperformed, on average, all other methods. The difference of the latter approach compared to \ks\ is noticeable mainly for data set {\ttfamily{TB5}}.
    \end{itemize}
    \item[\textit{(b)}] The average computing times required by \alg, \ks, and \cpx\ are comparable, and all considerably shorter than those for \hils. On average, \hils\ is roughly seven times slower than \alg. The difference in terms of computing efficiency is more evident as the size of the instances grows. Notice that for the largest data set {\ttfamily{TB5}}, \alg\ is, on average, more than ten times faster than \hils. 
\end{itemize}


\begin{table}
\centering
\caption{Benchmark instances: A comparison of the performance of \cpx, \ks, \hils, and \alg. \label{tab:AgainstLiterature}}
{
\setlength{\tabcolsep}{2.5pt}
\begin{tabular}{@{}lc rrr rrr rrr rrr@{}}
\hline
&  &    \multicolumn{3}{l}{\textbf{\cpx}}  &    \multicolumn{3}{l}{\textbf{\ks}} &    \multicolumn{3}{l}{\textbf{\hils}} &    \multicolumn{3}{l}{\textbf{\alg}}  \\ 
              \cline{3-5} \cline{6-8} \cline{9-11} \cline{12-14}
\textbf{Data}              & \textbf{\#} &      &  \multicolumn{1}{c}{$\bm{z}^{UB}$} &  \multicolumn{1}{c}{\textbf{Time}}   &        &  \multicolumn{1}{c}{$\bm{z}^{UB}$} &  \multicolumn{1}{c}{\textbf{Time}} &      &  \multicolumn{1}{c}{$\bm{z}^{UB}$} &  \multicolumn{1}{c}{\textbf{Time}} &       &  \multicolumn{1}{c}{$\bm{z}^{UB}$} &  \multicolumn{1}{c}{\textbf{Time}} \\ 
\textbf{Set}              & \textbf{Inst.} & \multicolumn{1}{c}{\textbf{\#Best}}     & \multicolumn{1}{c}{\textbf{Gap (\%)}} &  \multicolumn{1}{c}{\textbf{(sec.)}}    & \multicolumn{1}{c}{\textbf{\#Best}}     & \multicolumn{1}{c}{\textbf{Gap (\%)}} &  \multicolumn{1}{c}{\textbf{(sec.)}} & \multicolumn{1}{c}{\textbf{\#Best}}     & \multicolumn{1}{c}{\textbf{Gap (\%)}} &  \multicolumn{1}{c}{\textbf{(sec.)}} & \multicolumn{1}{c}{\textbf{\#Best}}     & \multicolumn{1}{c}{\textbf{Gap (\%)}} &  \multicolumn{1}{c}{\textbf{(sec.)}}\\
\hline 
{\ttfamily{OR4}} & 12 & 12    & 0.00 & 918   & 12    & 0.00 & 146   & 12    & 0.00 & 922   & 12    & 0.00 & 369 \\  
\hline
{\ttfamily{TB1}} & 20 & 19    & 0.00 & 1,031  & 17    & 0.01 & 1,330  & 17    & 0.01 & 3,593  & 15    & 0.03 & 1,040 \\
{\ttfamily{TB2}} & 20 & 20    & 0.00 & 383   & 20    & 0.00 & 236   & 20    & 0.00 & 663   & 20    & 0.00 & 329 \\
{\ttfamily{TB3}} & 20 &  15    & 0.04 & 2,062  & 16    & 0.02 & 2,162  & 14    & 0.03 & 10,819 & 14    & 0.04 & 2,312 \\
{\ttfamily{TB4}} & 20 & 7     & 0.17 & 2,947  & 11    & 0.05 & 2,855  & 7     & 0.09 & 22,467 & 15    & 0.04 & 2,950 \\
{\ttfamily{TB5}} & 20 & 4     & 0.45 & 3,595  & 8     & 0.11 & 3,602  & 4     & 0.17 & 36,606 & 14    & 0.03 & 3,550 \\ 
\cline{1-2} \cline{3-5} \cline{6-8} \cline{9-11} \cline{12-14}
\textbf{Avg. (Tot.)} & \textbf{(112)}  & \textbf{(77)} & \textbf{0.12} & \textbf{1,887} & \textbf{(84)} & \textbf{0.03} & \textbf{1,834} & \textbf{(74)} & \textbf{0.05} & \textbf{13,340} & \textbf{(90)} & \textbf{0.02} & \textbf{1,858}  \\ \hline
\end{tabular}}
\end{table}

To investigate the impact of the methodological contributions introduced in the present paper, we display in Table~\ref{tab:PaKSvsKS14} the value of some key statistics that explain the different performance of \alg\ compared to \ks\ on the largest benchmark instances. Each column from 3 to 7 shows a ratio (in \%) between the value that a given metric took for \alg\ and for \ks. We considered the following statistics. The time spent before entering the optimization step (see Section~\ref{sec:PKS} for \alg) is denoted by \textit{Time Pre.}. Statistic \textit{Init. $|K_y|$} indicates the number of $y$ variables composing the initial kernel $K_y$. Concerning the buckets, the number of buckets generated is denoted by \textit{\#$B^r$}, whereas \textit{Avg. $|B^r_y|$} and \textit{Avg. $|B^r_x|$} are the average number of $y$ and $x$ variables that compose the buckets, respectively. Finally, statistic \textit{$y$ remov.} denotes the number of $y$ variables removed during the initialization phase of \alg. Recall that no $y$ variable is removed in \ks. From Table~\ref{tab:PaKSvsKS14} we can gain the following main insights:

\begin{itemize}
    \item[\textit{(i)}] Not surprisingly, \alg\ spends, on average, 41\% more time than \ks\ to carry out all steps that are preliminary to the optimization step;
    \item[\textit{(ii)}] The number of $y$ variables composing the initial kernel $K_y$ for \alg\ is, on average, 30\% more than that for \ks. This, along with the following observation, suggests a better capacity of \alg\ than \ks\ to identify the promising variables;
    \item[\textit{(iii)}] \alg\ generates, on average, half the number of buckets created by \ks. The size of the buckets is, on average, smaller for \alg,  both in terms of the number of $y$ and $x$ variables. This leads to smaller and, frequently, easier to solve restricted BIPs;
    \item[\textit{(iv)}] \alg\ discards, on average, 61\% of the $y$ variables. Notice that when one of these variables is removed (i.e., in all restricted BIPs we fix to zero a given $y$ variable), all of the associated $x$ variables are discarded due to constraints \eqref{SSCFL-demand}.
\end{itemize}

\begin{table}
\centering
\caption{Benchmark instances: Evaluating the impact of Phase 1. \label{tab:PaKSvsKS14}}
{\begin{tabular}{@{\extracolsep{4pt}}lc rrrrr r}
\hline 
&  &    \multicolumn{5}{l}{\textbf{$\bm{\Delta}=$ \alg / \ks}}  &    \multicolumn{1}{l}{\textbf{\alg}}  \\ 
             \cline{3-7} \cline{8-8}
 \textbf{Data}              & \textbf{\#} &  \multicolumn{1}{c}{\textbf{Time}} &  \multicolumn{1}{c}{\textbf{Init.}}   &       &  \multicolumn{1}{c}{\textbf{Avg.}} &  \multicolumn{1}{c}{\textbf{Avg.}} &   \multicolumn{1}{c}{$\bm{y}$}   \\ 
 \textbf{Set}              & \textbf{Inst.} & \multicolumn{1}{c}{\textbf{Pre.}}     & \multicolumn{1}{c}{$\bm{|K_y|}$} &  \multicolumn{1}{c}{$\bm{\#B^r}$}    &  \multicolumn{1}{c}{$\bm{|B^r_y|}$}     &  \multicolumn{1}{c}{$\bm{|B^r_x|}$} &  \multicolumn{1}{c}{\textbf{remov.}} \\
 \hline 
{\ttfamily{OR4}} & 12 & 244\% & 100\% & 37\%  & 46\%  & 40\%  & 72\% \\ 
 \hline
 {\ttfamily{TB1}} & 20 & 195\% & 146\% & 48\%  & 55\%  & 50\%  & 69\% \\
 {\ttfamily{TB2}} & 20 & 116\% & 115\% & 125\% & 71\%  & 60\%  & 5\% \\
 {\ttfamily{TB3}} & 20 & 123\% & 139\% & 38\%  & 54\%  & 47\%  & 73\% \\
 {\ttfamily{TB4}} & 20 & 109\% & 137\% & 37\%  & 45\%  & 38\%  & 74\% \\
 {\ttfamily{TB5}} & 20 & 99\%  & 133\% & 30\%  & 52\%  & 45\%  & 77\% \\
 \cline{1-2} \cline{3-7} \cline{8-8} 
 \textbf{Avg. (Tot.)} & \textbf{(112)}  & \textbf{141\%} & \textbf{130\%} & \textbf{54\%}  & \textbf{55\%}  & \textbf{47\%}  & \textbf{61\%} \\ \hline
\end{tabular}}
\end{table}

\subsection{Performance evaluation on the new large-scale instances} \label{sec:NewInstances}

In this section, we assess the performance of \alg\ on the new large-scale instances. To validate the performance, we employ the same set of statistics previously introduced, with the following differences. Each new instance is solved with \cpx, \ks, and \alg , and $z^{UB}$ is the value of the best solution produced across all of them. As these data sets are novel for the \prob, and for the majority of them \cpx\ did not prove the optimality of the solution found, we add statistic \textit{$z^{LB}$ Gap (\%)}. This statistic is computed as $100 \times \frac{z^H - z^{LB}}{z^{LB}}$, where $z^{LB}$ is the value of the best lower bound produced by \cpx\ before reaching the time limit or terminating due to an out of memory error. Recall that, as explained in Section~\ref{sec:Testing}, for the latter reason computing times of \cpx\ are not reported in the following analysis. The main rationale for reporting statistic \textit{$z^{LB}$ Gap (\%)} is to identify the instances where \cpx\ found a near-optimal, or even optimal, solution and to determine how far the value of the best solution found by each method can be from the optimum.

\begin{figure}[h]
     \begin{minipage}[t]{0.5\textwidth}
\centering
\begin{tikzpicture}[scale=0.9]
\begin{axis}[
  ylabel={Time (sec.)},    
  xlabel={Best solutions found},
    scatter/classes={
        a1={mark=square*,blue,fill opacity=0.3},
        a2={mark=square*,magenta,fill opacity=0.3},
        b1={mark=triangle*,blue,fill opacity=0.3},
        b2={mark=triangle*,magenta,fill opacity=0.3},
        c1={mark=*,blue,fill opacity=0.3},
        c2={mark=*,magenta,fill opacity=0.3},
        d1={mark=x,blue,fill opacity=0.3},
        d2={mark=x,magenta,fill opacity=0.3}%
    },
] 
\addplot[
        scatter, 
        only marks,
        scatter src=explicit symbolic,
    ]
    table[meta=label] {
        x       y       label   annotation
        209      4959    a1      \ks    
        254      5050    a2      \alg
    };
    \addplot[
        scatter, 
        only marks,
        scatter src=explicit symbolic,
    ]
    table[meta=label] {
        x       y       label   annotation
        31      7168    b1      \ks    
        72      6772    b2      \alg
        56      6118    c1      \ks
        68      5934    c2      \alg
        122      2245    d1      \ks
        114      2948    d2      \alg
    };
\end{axis}
\matrix [draw, above left] at (6.5,4.1) {
  \node[blue,fill opacity=0.9,font=\tiny] {\ks}; \\
    \node[magenta,fill opacity=0.9,font=\tiny] { \alg }; \\
};
     \end{tikzpicture}
     \subcaption{\#Best vs. Time. ($\square$=All inst.; $\triangle$={\ttfamily{TB-A}}; $\medcirc$={\ttfamily{TB-B}}; $\times$={\ttfamily{TB-C}} -– \cpx\ not shown due to the excessive out of
memory crashes.)}\label{fig:ScattNewA}
\end{minipage}\hfill
\begin{minipage}[t]{0.5\textwidth}
   \begin{tikzpicture}[scale=0.9]
      \centering
  \begin{axis}
    [
    boxplot/draw direction=x,
    boxplot={
    draw position={1/4 + floor(\plotnumofactualtype/3) + 1/4*mod(\plotnumofactualtype,3)},
    %
    box extend=0.25,
    },
    ylabel={},
    xlabel={$z^{UB}$ Gap (\%)},
    cycle list={{magenta},{blue},{cyan}},
    y=1.75cm,
    ytick={0,1,2,3},
    y tick label as interval,
    yticklabels={{\ttfamily{TB-C}}, {\ttfamily{TB-B}}, {\ttfamily{TB-A}}},
    xticklabel={\pgfmathparse{\tick}\pgfmathprintnumber{\pgfmathresult}\%} 
    ]
    \addplot+[
    fill,fill opacity=0.3,
    boxplot prepared={
      median=0.00,
      upper quartile=0.00,
      lower quartile=0.00,
      upper whisker=0.00, 
      lower whisker=0.00
    },
    ] coordinates {};
    \addplot+[
    fill,fill opacity=0.3,
    boxplot prepared={
      median=0.00,
      upper quartile=0.00,
      lower quartile=0.00,
      upper whisker=0.00, 
      lower whisker=0.00
    },
    ] coordinates {};
    \addplot+[
    fill,fill opacity=0.3,
    boxplot prepared={
      median=0.00,
      upper quartile=0.00,
      lower quartile=0.00,
      upper whisker=0.00, 
      lower whisker=0.00
    },
    ] coordinates {};
    \addplot+[
    fill,fill opacity=0.3,
    boxplot prepared={
      median=0.00,
      upper quartile=0.07,
      lower quartile=0.00,
      upper whisker=0.15, 
      lower whisker=0.00
    },
    ] coordinates {};
        \addplot+[
    fill,fill opacity=0.3,
    boxplot prepared={
      median=0.00,
      upper quartile=0.15,
      lower quartile=0.00,
      upper whisker=0.34, 
      lower whisker=0.00
    },
    ] coordinates {};
        \addplot+[
    fill,fill opacity=0.3,
    boxplot prepared={
      median=0.77,
      upper quartile=1.49,
      lower quartile=0.06,
      upper whisker=3.63, 
      lower whisker=0.00
    },
    ] coordinates {};
        \addplot+[
    fill,fill opacity=0.3,
    boxplot prepared={
      median=0.00,
      upper quartile=0.08,
      lower quartile=0.00,
      upper whisker=0.18, 
      lower whisker=0.00
    },
    ] coordinates {};
        \addplot+[
    fill,fill opacity=0.3,
    boxplot prepared={
      median=0.13,
      upper quartile=0.55,
      lower quartile=0.00,
      upper whisker=1.29, 
      lower whisker=0.00
    },
    ] coordinates {};
        \addplot+[
    fill,fill opacity=0.3,
    boxplot prepared={
      median=2.37,
      upper quartile=4.18,
      lower quartile=0.71,
      upper whisker=7.55, 
      lower whisker=0.00
    },
    ] coordinates {};
    \end{axis}
    \matrix [draw, above left] at (6.5,0.5) {
  \node[cyan,fill opacity=0.9,font=\tiny] { \cpx}; \\
  \node[blue,fill opacity=0.9,font=\tiny] {\ks}; \\
    \node[magenta,fill opacity=0.9,font=\tiny] { \alg }; \\
};
    \end{tikzpicture}
    \subcaption{Box-plots of $z^{UB}$ Gap. (Outliers not shown.) }\label{fig:ScattNewB}
\end{minipage}
  \caption{New large-scale instances (All Feas.): A comparison of the solution approaches.}\label{fig:ScattNew}
\end{figure}
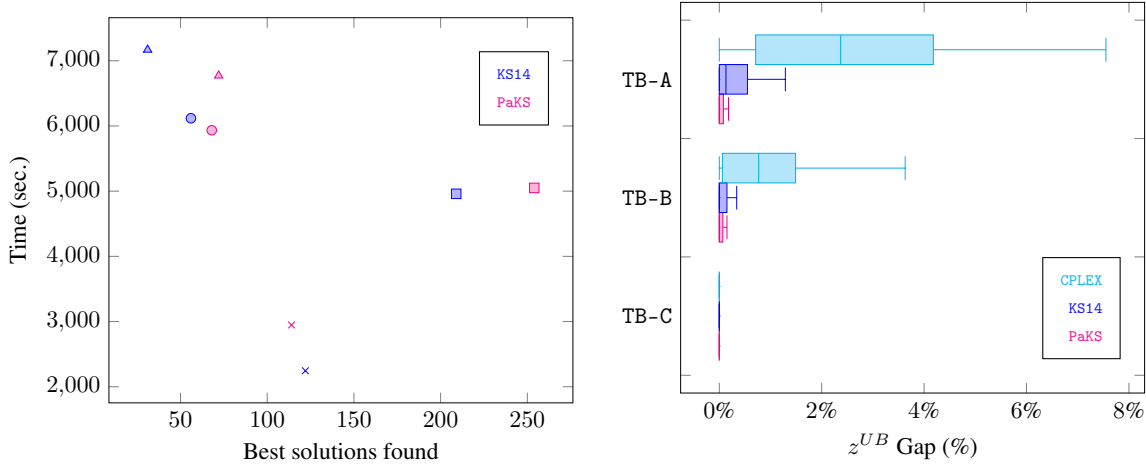

\begin{table}
\centering
\caption{New large-scale instances (All Feas.): A comparison of the performance of \cpx, \ks, and \alg. (Computing times of \cpx\ are not reported due to the excessive out of memory crashes.)  \label{tab:NewInstances}}
{
\setlength{\tabcolsep}{2.5pt}
\begin{tabular}{@{}lc rrr rrrr rrrr }
\hline 
&  &    \multicolumn{3}{l}{\textbf{\cpx}}  &    \multicolumn{4}{l}{\textbf{\ks}} &      \multicolumn{4}{l}{\textbf{\alg}}  \\ 
              \cline{3-5} \cline{6-9} \cline{10-13} 
\textbf{Data}              & \textbf{\#} &      &  \multicolumn{1}{c}{$\bm{z}^{UB}$} &  \multicolumn{1}{c}{$\bm{z}^{LB}$}   &        &  \multicolumn{1}{c}{$\bm{z}^{UB}$} &  \multicolumn{1}{c}{$\bm{z}^{LB}$}   &  \multicolumn{1}{c}{\textbf{Time}} &      &  \multicolumn{1}{c}{$\bm{z}^{UB}$} &  \multicolumn{1}{c}{$\bm{z}^{LB}$}   &  \multicolumn{1}{c}{\textbf{Time}}  \\ 
\textbf{Set}              & \textbf{Inst.} & \multicolumn{1}{c}{\textbf{\#Best}}     & \multicolumn{1}{c}{\textbf{Gap (\%)}} & \multicolumn{1}{c}{\textbf{Gap (\%)}}  & \multicolumn{1}{c}{\textbf{\#Best}}     & \multicolumn{1}{c}{\textbf{Gap (\%)}} & \multicolumn{1}{c}{\textbf{Gap (\%)}} &  \multicolumn{1}{c}{\textbf{(sec.)}} & \multicolumn{1}{c}{\textbf{\#Best}}    & \multicolumn{1}{c}{\textbf{Gap (\%)}} & \multicolumn{1}{c}{\textbf{Gap (\%)}} &  \multicolumn{1}{c}{\textbf{(sec.)}} \\
\hline 
{\ttfamily{TB-A1}} & 25    & 1     & 5.27 & 7.02 & 8     & 2.53 & 4.25 & 6,978  & 20    & 0.14 & 1.79 & 7,090 \\
{\ttfamily{TB-A2}} & 22    & 2     & 0.63 & 3.99 & 8     & 0.32 & 3.67 & 7,288  & 13    & 0.25 & 3.59 & 7,207 \\
{\ttfamily{TB-A3}} & 20    & 1     & 3.44 & 5.56 & 5     & 3.15 & 5.28 & 7,055  & 15    & 0.12 & 2.16 & 6,686 \\
{\ttfamily{TB-A4}} & 18    & 0     & 4.33 & 6.83 & 3     & 0.88 & 3.38 & 7,296  & 15    & 0.02 & 2.46 & 7,222 \\
{\ttfamily{TB-A5}} & 19    & 3     & 2.96 & 8.11 & 7     & 0.99 & 6.14 & 7,277  & 9     & 0.67 & 5.76 & 5,515 \\
\cline{1-2} \cline{3-5} \cline{6-9} \cline{10-13}
\textbf{Avg. (Tot.)} & \textbf{(104)} & \textbf{(7)} & \textbf{3.35} & \textbf{6.26} & \textbf{(31)} & \textbf{1.62} & \textbf{4.52} & \textbf{7,168} & \textbf{(72)} & \textbf{0.24} & \textbf{3.08} & \textbf{6,772} \\ 
\hline
{\ttfamily{TB-B1}} & 23    & 2     & 0.96 & 3.05 & 11    & 13.65 & 16.39 & 5,975  & 15    & 0.12 & 2.20 & 5,766 \\
{\ttfamily{TB-B2}} & 21    & 12    & 0.30 & 1.33 & 9     & 0.61 & 1.65 & 5,635  & 16    & 0.06 & 1.09 & 5,620 \\
{\ttfamily{TB-B3}} & 25    & 4     & 1.48 & 3.24 & 18    & 0.09 & 1.82 & 5,618  & 18    & 0.21 & 1.94 & 5,641 \\
{\ttfamily{TB-B4}} & 19    & 5     & 9.33 & 14.51 & 9     & 3.43 & 7.94 & 6,836  & 11    & 1.46 & 5.76 & 5,875 \\
{\ttfamily{TB-B5}} & 16    & 1     & 1.98 & 4.01 & 9     & 0.11 & 2.12 & 6,882  & 8     & 0.20 & 2.21 & 7,114 \\
\cline{1-2} \cline{3-5} \cline{6-9} \cline{10-13}
\textbf{Avg. (Tot.)} & \textbf{(104)} & \textbf{(24)} & \textbf{2.64} & \textbf{4.99} & \textbf{(56)} & \textbf{3.81} & \textbf{6.17} & \textbf{6,118} & \textbf{(68)} & \textbf{0.39} & \textbf{2.56} & \textbf{5,934} \\
\hline
{\ttfamily{TB-C1}} & 25    & 23    & 0.07 & 0.41 & 24    & 0.06 & 0.40 & 2,377  & 22    & 0.01 & 0.35 & 2,636 \\
{\ttfamily{TB-C2}} & 26    & 26    & 0.00 & 0.00 & 26    & 0.00 & 0.00 & 614   & 22    & 0.13 & 0.13 & 2,104 \\
{\ttfamily{TB-C3}} & 25    & 22    & 0.30 & 1.42 & 24    & 0.00 & 1.10 & 2,948  & 22    & 0.39 & 1.54 & 3,505 \\
\hdashline
{\ttfamily{TB-C4}} & 28    & 25    & 3.86 & 5.81 & 23    & 25.00 & 29.04 & 3,046  & 25    & 13.13 & 16.85 & 3,301 \\
(no 3 outl.) & 25    & 22    & 4.32 & 4.91 & 23    & 13.04 & 13.99 & 2,546  & 25    & 0.00 & 0.41 & 2,829 \\
\hdashline
{\ttfamily{TB-C5}} & 25    & 25    & 0.00 & 0.01 & 25    & 0.00 & 0.01 & 2,208  & 23    & 0.01 & 0.01 & 3,187 \\
\cline{1-2} \cline{3-5} \cline{6-9} \cline{10-13}
\textbf{Avg. (Tot.)} & \textbf{(129)} & \textbf{(121)} & \textbf{0.91} & \textbf{1.62} & \textbf{(122)} & \textbf{5.44} & \textbf{6.60} & \textbf{2,245} & \textbf{(114)} & \textbf{2.95} & \textbf{4.05} & \textbf{2,948} \\ 
\hdashline
\textbf{(no 3 outl.)} & \textbf{(126)} & \textbf{(118)} & \textbf{0.93} & \textbf{1.34} & \textbf{(122)} & \textbf{2.60} & \textbf{3.08} & \textbf{2,126} & \textbf{(114)} & \textbf{0.11} & \textbf{0.48} & \textbf{2,846} \\ 
\hline
\hline
\textbf{Avg. (Tot.)} & \textbf{(337)} & \textbf{(152)} & \textbf{2.20} & \textbf{4.09} & \textbf{(209)} & \textbf{3.76} & \textbf{5.83} & \textbf{4,959} & \textbf{(254)} & \textbf{1.32} & \textbf{3.29} & \textbf{5,050} \\ 
\hdashline
\textbf{(no 3 outl.)} & \textbf{(334)} & \textbf{(149)} & \textbf{2.22} & \textbf{4.01} & \textbf{(209)} & \textbf{2.67} & \textbf{4.49} & \textbf{4,939} & \textbf{(254)} & \textbf{0.23} & \textbf{1.94} & \textbf{5,030} \\ 
\hline
\end{tabular}}
\end{table}

Figure~\ref{fig:ScattNew} and Table~\ref{tab:NewInstances} summarize the results for the new large-scale instances. For readability, Figures~\ref{fig:ScattNewA} and \ref{fig:ScattNewB} show the results clustered according to the three main data sets {\ttfamily{TB-A}}, {\ttfamily{TB-B}}, and {\ttfamily{TB-C}}. Figure~\ref{fig:ScattNewA} also shows the results for all new instances (depicted as squares), whereas outliers are not displayed in Figure~\ref{fig:ScattNewB}. For each of the main data sets, the results are detailed in Table~\ref{tab:NewInstances} according to the size of the instances.

The main insights we can gain from Figure~\ref{fig:ScattNew} and Table~\ref{tab:NewInstances} are as follows:

\begin{itemize}
    \item[\textit{(a)}] \alg\ outperforms both \cpx\ and \ks\ in terms of quality of the solution found: 
    \begin{itemize}
        \item[\textit{(i)}] \alg\ produced high-quality solutions, with 254 out of the 337 best solutions, a value significantly higher than \ks\ (209) and \cpx\ (152);        
        \item[\textit{(ii)}] \alg\ achieved the smallest average value of \textit{$z^{UB}$ Gap (\%)} of 1.32\%, which is considerably better than both \cpx\ (2.20\%) and \ks\ (3.76\%);
        \item[\textit{(iii)}] \alg\ achieved the largest improvements, compared to the two other approaches, on data sets {\ttfamily{TB-A}} and {\ttfamily{TB-B}}. For these two data sets, the average values of \textit{$z^{UB}$ Gap (\%)} for \alg\ is 0.24\% and 0.39\%, respectively, which are remarkably smaller than the corresponding average values obtained by \cpx\ and \ks;
        \item[\textit{(iv)}] Regarding the instances of data set {\ttfamily{TB-C}} (see Table~\ref{tab:NewInstances}), \cpx\ outperformed, on average, both other methods. Nevertheless, this outcome is only caused by 3 out of the 129 instances in this data set where \alg\ produced solutions considerably worse than those found by \cpx\ (additional details can be found below). Without outliers, the three approaches perform equally well (see Figure~\ref{fig:ScattNewB}).
    \end{itemize}
    \item[\textit{(b)}] Computing times for \alg\ and \ks\ are, on average, comparable.
\end{itemize}

As previously mentioned, if we compare the performance of the three methods for all the instances of data set {\ttfamily{TB-C}}, \cpx\ achieved the smallest average value of \textit{$z^{UB}$ Gap (\%)}. We investigated this issue in more detail. First, it is worth observing that for several of the instances in this data set, \cpx\ found a proven optimal (or near-optimal) solution, except for the instances of data set {\ttfamily{TB-C4}} (observe the values of \textit{$z^{LB}$ Gap (\%)} in Table~\ref{tab:NewInstances}). Three of the latter instances turned out to be particularly challenging for \alg, which produced solutions that deteriorated significantly with respect to those found by \cpx. In addition to showing in Figure~\ref{fig:ScattNewB} the distribution of \textit{$z^{UB}$ Gap (\%)} without outliers, we report in Table~\ref{tab:NewInstances} the value of each statistic computed after removing the 3 outliers --rows with `\textit{(no 3 outl.)}'. Without those 3 outliers, the average values of  \textit{$z^{UB}$ Gap (\%)} achieved by \alg\ have sharp improvements: from 13.13\% to 0.00\% for the instances of data set {\ttfamily{TB-C4}}; from 2.95\% to 0.11\% for {\ttfamily{TB-C}}; and from 1.32\% to 0.23\% for the full set of new large-scale instances.

\section{Conclusions}\label{sec:Conclusions}

We introduced and computationally validated a two-phase matheuristic, called the Pattern-based Kernel Search (\alg), for the solution of the Single-Source Capacitated Facility Location Problem (\prob). \alg\ extends the conventional Kernel Search (KS) framework by combining the latter with a pattern recognition technique employed to gather information on the importance and interdependence among decision variables. An extensive computational study conducted on benchmark and new very large-scale instances was carried out. On the benchmark instances, \alg\ was validated against \cpx\ and two other heuristics from the literature, a standard KS re-implementation (\ks) and a Hybrid Iterated Local Search (\hils), the latter being the current state-of-the-art algorithm for solving the \prob. For these instances, comprising up to 1,000 locations and 1,000 customers, \alg\ outperformed, in terms of quality of the solutions produced, both \ks\ and \hils, as well as \cpx\ when run with a time limit. Regarding the computing times, \alg\ was consistently faster than \hils, and comparable to \ks\ and \cpx. On the new very large-scale instances, consisting of up to 2,000 locations and 2,000 customers, \alg\ outperformed both \ks\ and \cpx\ in terms of quality of the solutions found. Further, it outperformed both methods also in terms of robustness, encountering remarkably fewer problems than the other two approaches --such as out of memory issues, or termination without any feasible integer solution. Computing times of \alg\ were comparable to those of \ks.

The present research illustrates the potential advantages of combining optimization algorithms with approaches from data science, which can make the resulting solution method data-driven and adaptive to some pivotal characteristics of the input data. Potential research avenues include the application of similar combined solution approaches to other classes of combinatorial optimization problems.

\textbf{\textit{Acknowledgements}}
Research visits of the collaborating parties were financially supported by the Karlsruhe House of Young Scientists (KHYS). The authors acknowledge support by the state of Baden-Württemberg through bwHPC.\\ The research of G. Guastaroba and M.G. Speranza has been partially supported by the following projects funded by the European Union (EU) and the Italian Ministry for Universities and Research (MUR): (1) National Recovery and Resilience Plan (NRRP), within project “Sustainable Mobility Center (MOST)” (2022–2026), CUP: D83C22000690001, Spoke N. 7, ``CCAM, Connected networks and Smart Infrastructures”; (2) Project PRIN 2022 - Next Generation EU ``Time-dependent optimization for sustainable transportation", Code: 20223MHHA8, CUP: D53D23005590006.

\bibliographystyle{informs2014} 
\bibliography{Literature} 






\end{document}